\documentclass{amsart}

\usepackage{amsfonts}
\usepackage{amsmath}
\usepackage{amsthm}
\usepackage{amssymb}
\usepackage[english]{babel}
\usepackage{graphics}
\usepackage{latexsym}
\usepackage{longtable} \usepackage{mathrsfs}
\usepackage{graphicx}
\usepackage{wrapfig} 
\usepackage[pdftex,colorlinks=true]{hyperref}

\newtheorem{theorem}{Theorem}
\newtheorem{corollary}[theorem]{Corollary}
\newtheorem{lemma}[theorem]{Lemma}
\newtheorem{proposition}[theorem]{Proposition}

\newtheorem{claim}[theorem]{Claim}
\newtheorem{example}[theorem]{Example}

\theoremstyle{definition}
\newtheorem{definition}[theorem]{Definition}
\newtheorem{remark}[theorem]{Remark}

\newcommand{\msL}{\mathscr{L}}

\newcommand{\mH}{\mathcal{H}}

\newcommand{\bR}{\bold{R}}

\renewcommand{\S}{\mathcal{S}}

\newcommand{\A}{\textrm{A}}
\newcommand{\R}{\mathbb{R}}
\newcommand{\N}{\mathbb{N}}
\newcommand{\mS}{\mathbb{S}}
\newcommand{\mB}{\mathbb{B}}

\renewcommand{\P}{\textrm{P}}

\newcommand{\noi}{\noindent}
\newcommand{\ms}{\medskip}

\newcommand{\al}{\alpha}
\newcommand{\be}{\beta}
\newcommand{\ga}{\gamma}
\newcommand{\de}{\delta}
\newcommand{\De}{\Delta}
\newcommand{\e}{\varepsilon}
\newcommand{\si}{\sigma}
\newcommand{\la}{\lambda}
\newcommand{\ka}{\kappa}
\newcommand{\Om}{\Omega}

\newcommand{\lharpoonup}{-\!\!\!-\!\!\!\!\rightharpoonup}
\newcommand{\weakstar}{\overset{ \phantom{a}_* \ }{\lharpoonup}}
\newcommand{\larrow}{\longrightarrow}
\newcommand{\ot}{\otimes}

\newcommand{\ri}{\rightarrow}

\newcommand{\p}{\partial}
\newcommand{\sub}{\subseteq}
\newcommand{\set}{\setminus}
\newcommand{\by}{\times}
\newcommand{\rk}{\textrm{rk}}

\newcommand{\tr}{\textrm{tr}}
\newcommand{\sgn}{\textrm{sgn}}

\newcommand{\diam}{\textrm{diam}}
\newcommand{\dist}{\textrm{dist}}

\newcommand{\Div}{\textrm{Div}}
\newcommand{\Curl}{\textrm{Curl}}
\newcommand{\co}{\overline{\textrm{co}}}
\newcommand{\inter}{\textrm{int}}
\newcommand{\spn}{\textrm{span}}

\newcommand{\bt}{\begin{theorem}}\newcommand{\et}{\end{theorem}}
\newcommand{\bd}{\begin{definition}}\newcommand{\ed}{\end{definition}}
\newcommand{\bl}{\begin{lemma}}\newcommand{\el}{\end{lemma}}
\newcommand{\beq}{\begin{equation}}\newcommand{\eeq}{\end{equation}}
\newcommand{\bc}{\begin{claim}}\newcommand{\ec}{\end{claim}}
\newcommand{\bex}{\begin{example}}\newcommand{\eex}{\end{example}}
\newcommand{\bcor}{\begin{corollary}}\newcommand{\ecor}{\end{corollary}}
\newcommand{\bp}{\begin{proof}}\newcommand{\ep}{\end{proof}}

\newcommand{\BPCOR}{\medskip \noindent \textbf{Proof of Corollary} }
\newcommand{\BPP}{\medskip \noindent \textbf{Proof of Proposition} }
\newcommand{\BPT}{\medskip \noindent \textbf{Proof of Theorem} }

\numberwithin{equation}{section}
\numberwithin{theorem}{section}

\begin{document}

\title{On the Structure of $\infty$-Harmonic maps}

\author{\textsl{Nicholas Katzourakis}}
\address{Department of Mathematics and Statistics, University of Reading, Whiteknights, PO Box 220, Reading RG6 6AX, UK {\tiny and} BCAM - Basque Center for Applied Mathematics, Mazarredo 14, E48009, Bilbao, Spain}
\email{n.katzourakis@reading.ac.uk}

\subjclass[2010]{Primary 35J47, 35J62, 53C24; Secondary 49J99}

\date{}


\keywords{$\infty$-Laplacian, Aronsson equation, Rigidity Theory, Quasiconformal maps, Calculus of Variations in $L^\infty$, Optimal Lipschitz Extensions.}

\begin{abstract} Let $H \in C^2(\mathbb{R}^{N \times n})$, $H\geq 0$. The PDE system
\[  \label{1}
\A_\infty u \, :=\, \Big(H_P \ot H_P + H [H_P]^\bot H_{PP} \Big)(Du) : D^2 u\, = \, 0 \tag{1}
\]
arises as the ``Euler-Lagrange PDE" of vectorial variational problems for the functional $E_{\infty}(u,\Om) = \| H(Du) \|_{L^\infty(\Om)}$ defined on maps $u : \Omega \subseteq \mathbb{R}^n \longrightarrow \mathbb{R}^N$. \eqref{1} first appeared in the author's recent work \cite{K3}. The scalar case though has a long history initiated by Aronsson in \cite{A1}. Herein we study the solutions of \eqref{1} with emphasis on the case of  $n=2\leq N$ with $H$ the Euclidean norm on $\mathbb{R}^{N \times n}$, which we call the  ``$\infty$-Laplacian". By establishing a rigidity theorem for rank-one maps of independent interest, we analyse a phenomenon of separation of the solutions to phases with qualitatively different behaviour. As a corollary, we extend  to $N \geq 2$ the Aronsson-Evans-Yu theorem regarding non-existence of zeros of $|Du|$ and prove a Maximum Principle. We further characterise all $H$ for which \eqref{1} is elliptic and also study the initial value problem for the ODE system arising for $n=1$ but with $H(\cdot,u,u')$ depending on all the arguments.

\end{abstract}

\maketitle

\section{Introduction} \label{section1}

Let $H \in C^2(\R^{N \by n})$ be a nonegative function which we call Hamiltonian. In this paper we study the classical solutions $u :\Om \sub \R^n \larrow \R^N$ of the PDE system
\beq \label{1.7}
\A_\infty u\ :=\ \Big(H_\P \ot H_\P + H[H_\P]^\bot H_{\P \P}\Big)(Du):D^2u\ = \ 0.
\eeq
Here $[H_\P(P)]^\bot$ denotes the orthogonal projection on the nullspace of $H_\P(P)^\top  : \R^N \larrow \R^n$ and $H_P$ is the derivative matrix (for details see the Preliminaries \ref{Preliminaries}). The system \eqref{1.7} arises as a sort of Euler-Lagrange PDE of vectorial variational problems in $L^\infty$ for the functional
 \beq \label{1.8}
E_{\infty}(u,\Om)\, := \, \big\| H(Du)\|_{L^\infty(\Om)}.
 \eeq
Calculus of Variations in $L^\infty$ is very important for applications, since minimisation of the maximum value leads to more realistic models when compared to the more classical case of integral functionals in which case we minimise the average. \eqref{1.7} is a quasilinear 2nd order system in non-divergence form which was first formally derived  by the author in the recent work \cite{K3} as the limit of Euler-Lagrange equations of the functionals $\int_\Om \big(H(Du)\big)^p$  as $p\ri \infty$. Herein particular emphasis will be given on the 2D case for $n=2\leq N$ with $H(P)=\frac{1}{2}|P|^2$, where $|\cdot|$ is the Euclidean norm on $\R^{N \by n}$. In this case \eqref{1.7} simplifies to
\beq  \label{1.1}
\De_\infty u \ :=\ \Big(Du \ot Du  + |Du|^2 [Du]^\bot \! \ot I \Big) : D^2 u\ = \ 0. 
\eeq
We call \eqref{1.1} the ``$\infty$-Laplacian" and its solutions $\infty$-Harmonic maps. The name stems from its derivation which we now recall. After expansion and normalisation of the $p$-Laplace system $\De_p u= \Div\big(|Du|^{p-2}Du\big)=0$, we have
\beq  \label{1.2}
Du \ot Du : D^2 u \, +\, \frac{|Du|^2}{p-2}\De u\ = \ 0. 
\eeq
Let $[Du]^\top$ and $[Du]^\bot$ denote the orthogonal projections on the range of $Du$ and the nullspace of $Du^\top$ respectively. Since $[Du]^\top + [Du]^\bot =I$, by expanding $\De u$ with respect to these projections, we get
\beq  \label{1.3}
Du \ot Du : D^2 u \, +\, \frac{|Du|^2}{p-2}[Du]^\top \De u\ = \ -\frac{|Du|^2}{p-2}[Du]^\bot \De u. 
\eeq
By orthogonality, right and left hand side of \eqref{1.3} are normal to each other. Hence, they both vanish and \eqref{1.3} actually decouples to 2 systems. By renormalising the right hand side of \eqref{1.3}  and rearranging, we get
\beq  \label{1.4}
Du \ot Du : D^2 u  + |Du|^2[Du]^\bot \De u\ = \ -\frac{|Du|^2}{p-2}[Du]^\top \De u. 
\eeq
As $p \ri \infty$, \eqref{1.4} formally leads to \eqref{1.1}. In the case of \eqref{1.1} the projection $[Du]^\bot$ coincides with the projection on the geometric normal space of the image of the solution.  When $n=1$, the system simplifies to
\begin{align} \label{1.6a}
\De_\infty u \ & =\ (u' \ot u') u''\, +\, |u'|^2\Big(I - \frac{u' \ot u'}{|u'|^2} \Big)u''\ =\ |u'|^2u''.
\end{align}
In particular, it follows that $\infty$-Harmonic curves are \emph{affine} and no interesting phenomena arise.

When $N=1$, the normal coefficient $|Du|^2[Du]^\bot$ vanishes identically and the same holds when $u$ is submersion. The single $\infty$-Laplacian PDE $D_i u D_j u D^2_{ij}u=0$ and the related scalar $L^\infty$-variational problems have a long history. $\De_\infty$ was first derived and studied by Aronsson in the '60s in \cite{A3, A4} and has been extensively studied ever since (see for example Crandall \cite{C}, Barron-Evans-Jensen \cite{BEJ} and references therein). A major difficulty in its study is its degeneracy and the emergence of singular solutions (see e.g.\ \cite{A6, A7, K1}). In the last 25 years the single PDE has been studied in the context of Viscosity Solutions.

A further difficulty of the vectorial case which is not present in the scalar case is that  \emph{\eqref{1.7} has discontinuous coefficients even when the operator $A_\infty$ is applied to $C^\infty$ maps} which are solutions. As an example consider 
\beq \label{1.13}
u(x,y)\, :=\, e^{ix}-e^{iy} \ ,\ \ \ u \ :\  \R^2 \larrow \R^2. 
\eeq
In \cite{K3} we showed that \eqref{1.13} is a smooth solution of the $\infty$-Laplacian near the origin. However, the coefficient $|Du|^2[Du]^\bot$ of \eqref{1.1} is discontinuous. The problem is that the projection $[Du]^\bot$ ``jumps" when the dimension of the image changes. Indeed, for \eqref{1.13} we have $\rk(Du)=2$ off the diagonal $\{x=y\}$, while $\rk(Du)=1$ otherwise. Hence, the domain of \eqref{1.13} splits to 3 components, the ``2D phase $\Om_2$", whereon $u$ is essentially 2D, the ``interface $\S$" where the coefficients of $\De_\infty$ become discontinuous and  the ``1D phase $\Om_1$", whereon $u$ is essentially 1D (and in this case is empty). Much more intricate examples of smooth 2D $\infty$-Harmonic maps whose interfaces have triple junctions and corners are constructed in \cite{K6}. For any $K\in C^1(\R)$ with $\|K\|_{L^\infty(\R)} <\frac{\pi}{2}$,  the formula
\beq \label{2.2aa}
u(x,y)\ :=\ \int_y^x e^{iK(t)}dt
\eeq
defines a smooth $\infty$-Harmonic map whose phases are as shown in Figures 1(a), 1(b) below, when $K$ qualitatively behaves as shown in the Figures 2(a), 2(b) respectively.
\begin{align}
& \underset{\scriptstyle{\text{Figure 1(a). \hspace{100pt} Figure 1(b).\ \ \  }}}{ \includegraphics[scale=.18]{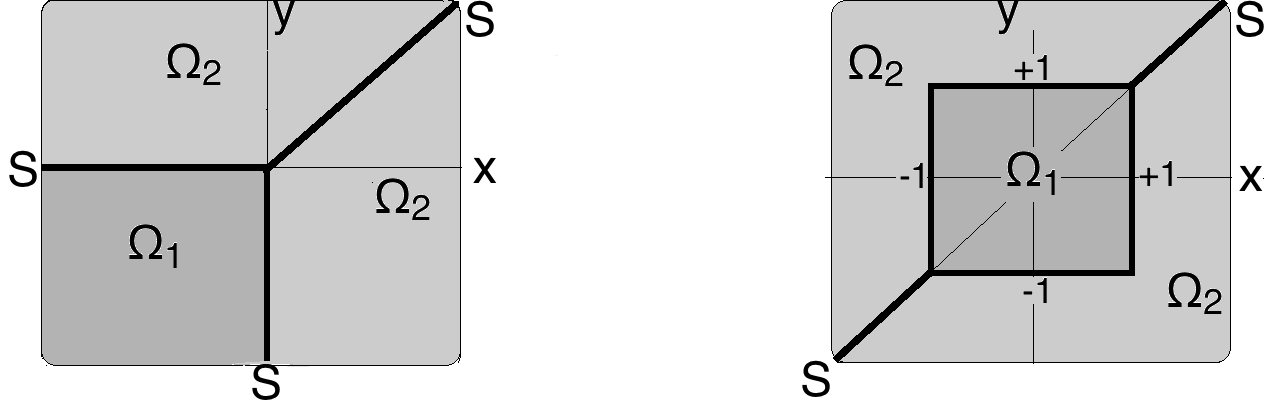}} \nonumber\\
& \underset{\scriptstyle{\text{Figure 2(a). \hspace{100pt} Figure 2(b).}}}{  \includegraphics[scale=.16]{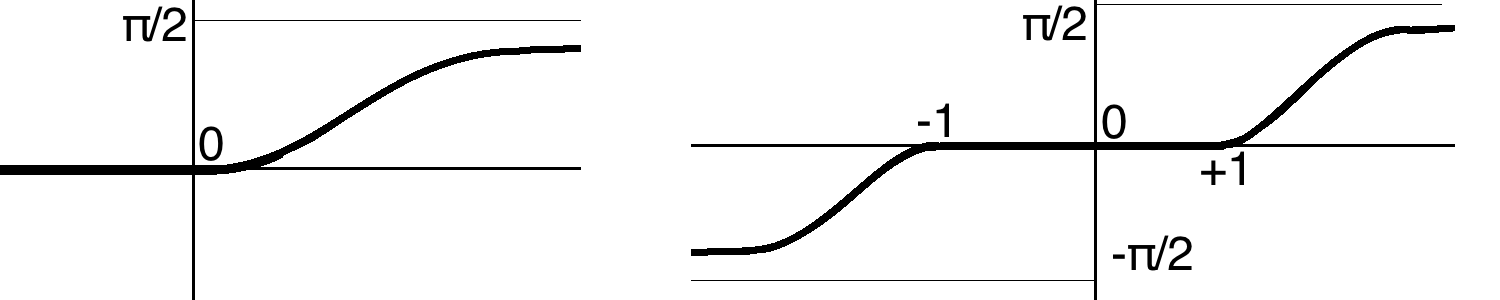} } \nonumber
\end{align}
Moreover, on $\Om_1$ \eqref{2.2aa} is given by a scalar $\infty$-Harmonic function times a constant vector, and on $\Om_2$ it is a solution of the vectorial Eikonal equation.

One of the principal results of this paper is that this \emph{phase separation} is a general phenomenon for smooth 2D $\infty$-Harmonic maps. On each phase the dimension of the tangent space is constant and these phases are separated by \emph{interfaces} whereon $[Du]^\bot$ becomes discontinuous. More precisely,  in Section \ref{section3}  we prove the next

\bt[Structure of 2D $\infty$-Harmonic maps] \label{th3} Let $u : \Om \sub \R^2 \larrow \R^N$ be an $\infty$-Harmonic map in $C^2(\Om)^N$, that is a solution to \eqref{1.1}. Let also $N\geq 2$. Then, there exists disjoint open sets $\Om_1$, $\Om_2 \sub \Om$ and a closed nowhere dense set $\S$ such that $\Om = \Om_1 \cup \S \cup \Om_2$ and:

\ms 

\noi (i) On $\Om_2$ we have $\rk(Du)=2$ and the map $u: \Om_2\larrow \R^N$ is an immersion and solution of the vectorial Eikonal equation:
\beq \label{3.2}
 |Du|^2\, = \, c^2\, > 0.
\eeq
The constant $c$ may vary on different connected components of $\Om_2$.
\ms

\noi (ii) On $\Om_1$ we have $\rk(Du)=1$ and the map $u: \Om_1\larrow \R^N$ is given by an essentially scalar $\infty$-Harmonic function $f  : \Om_1 \larrow \R$: 
\beq \label{3.3}
u \, = \, a + \xi f\ , \ \ \ \De_\infty f \, =\, 0,\ \ a\in \R^N, \ \xi \in \mS^{N-1}.
\eeq
The vectors $a,\xi$ may vary on different connected components of $\Om_1$.

\ms
\noi (iii) On $\S$, $|Du|^2$ is constant and also $\rk(Du)=1$. Moreover if $\S=\p \Om_1 \cap \p \Om_2$ (that is if both the 1D and 2D phases coexist) then  $u : \S\larrow \R^N$ is given by an essentially scalar solution of the Eikonal equation: 
\beq \label{3.4}
u \, = \, a + \xi f\ , \ \ \ |Df|^2 \, =\, c^2\, >0,\ \ a\in \R^N,\ \xi \in \mS^{N-1}.
\eeq
\et
\noi We note that this phase separation is a genuinely vectorial phenomenon, which does not arise when the rank is one. By employing Aronsson's result on the non-existence of zeros for the gradient of scalar $\infty$-Harmonic functions contained in \cite{A4}, we deduce the following consequence of Theorem \ref{th3}:
\begin{corollary}
[$\infty$-Harmonic maps have positive rank] \label{Cor3}  Let $u : \Om \sub \R^2 \larrow \R^N$ be an $\infty$-Harmonic map in $C^2(\Om)^N$. Then, either $|Du|>0$ on $\Om$ or $|Du|\equiv 0$ on $\Om$. Hence, non-constant $\infty$-Harmonic maps have positive rank.
\end{corollary}
\noi Corollary \ref{Cor3} is an extension to the vector case of the aforementioned theorem of Aronsson, which has been subsequently improved by Evans \cite{E} and Yu \cite{Y}. Hence, $\infty$-Harmonic maps have positive rank but generally non-constant rank. As a corollary,  in Section \ref{section3} we also establish a vectorial version of the Maximum Principle known as the \emph{Convex Hull Property}, valid for $n=N=2$:
\begin{corollary}[Convex Hull Property] 
\label{Cor2} 
Suppose that $u :\Om \sub \R^2 \larrow \R^2$ is an $\infty$-Harmonic map. Then, for all $\Om' \Subset \Om$, the image $u(\Om')$ is contained in the closed convex hull of the boundary values:
 \beq \label{3.17}
u(\Om') \ \sub \ \co\, \big(u(\p \Om')\big).
\eeq
\end{corollary}
\noi Since a convex set coincides with the intersection of half-spaces containing it, \eqref{3.17} is just an elegant formulation of the Maximum Principle for all 1D projections of $u$. It is well known in the context of Minimal Surfaces (see e.g.\ \cite{CM}, \cite{O}) and more generally in Calculus of Variations (see \cite{K2} and references therein). A topological consequence of Corollary \ref{Cor2} is

\begin{corollary}[Absence of interfaces] \label{Cor4} Suppose that $u :\Om \sub \R^2 \larrow \R^2$ is an $\infty$-Harmonic map. Then:

\ms

\noi (i) If $\Om_2 \Subset \Om$, then $\Om_2 = \emptyset$ and $\S = \emptyset$. Hence, either the set whereon $u$ is a local diffeomorphism has a common boundary portion with $\Om$ or it is empty and $u$ is everywhere essentially scalar without any interface $\S$.

\ms

\noi (ii) If $\Om\Subset \R^2$ and $u$ is essentially scalar near $\p \Om$, then there is no interface $\S$ inside $\Om$ and $u$ is essentially scalar throughout $\Om$.
\end{corollary}

The main analytical machinery required for the proof of Theorem \ref{th3} is developed in Section \ref{section2} and is a \emph{rigidity result for maps with 1D range} of independent interest. To begin with, consider a map $u : \Om \sub \R^n \larrow \R^N$ given as composition of a scalar function $f\in C^2(\Om)$ with a unit speed curve $\nu : \R \larrow \R^N$, that is $u = \nu \circ f$. Then, we have $Du = (\dot{\nu} \circ f)\ot Df$ and hence $u$ is a \emph{Rank-One map}, that is $\rk(Du)\leq 1$ on $\Om$.
\[
\underset{\text{Figure 3.}}{\includegraphics[scale=0.2]{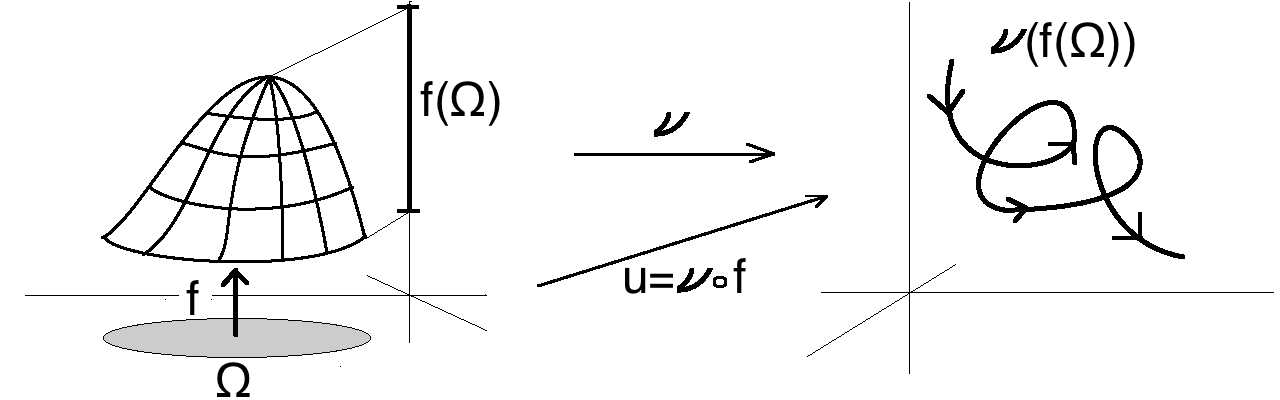}}
\]
Interestingly,  the class of Rank-One maps is  \emph{rigid} since a certain converse is true as well: all maps which satisfy $\rk(Du)\leq 1$ arise as compositions of unit speed curves with scalar functions. More precisely, 

 \bt[Rigidity of Rank-One maps] \label{th1} Suppose $\Om \sub \R^n$ is open and contractible and $u : \Om \sub \R^n \larrow \R^N$ is in $C^2(\Om)^N$. Then, the following are equivalent:

\ms 

\noi (i) $u$ is a Rank-One map, that is $\rk(Du)\leq 1$ on $\Om$ or equivalently there exist maps $\xi : \Om \larrow \R^N$ and $w : \Om \larrow \R^n$ with $w\in C^1(\Om)^n$ and $\xi \in C^1(\Om\set \{w=0\})^N$  such that $Du = \xi \ot w$.

\ms 

\noi (ii) There exists $f \in C^2(\Om)$, a partition $\{B_i\}_{i \in \N}$ of $\Om$ to Borel sets where each $B_i$ equals a connected open set with a boundary portion and Lipschitz curves $\{\nu^i \}_{i \in \N} \sub W^{1,\infty}_{loc}(\R)^N$ such that on each $B_i$ $u$ equals the composition of $\nu^i$ with $f$:
 \beq \label{2.1}
 u\ = \ \nu^i \circ f\ , \ \ \text{ on }B_i \sub \Om.
  \eeq 
Moreover, $|\dot{\nu}^i| \equiv 1$ on $f(B_i)$, $\dot{\nu}^i \equiv 0$ on $\R \set f(B_i)$ and there exists $\ddot{\nu}^i$ on $f(B_i)$, interpreted as $1$-sided on $\p f(B_i)$, if any. Also, 
 \beq \label{2.2}
  Du \ =\ (\dot{\nu}^i \circ f ) \ot Df \ , \ \ \text{ on }B_i \sub \Om,
 \eeq
and the image $u(\Om)$ is an $1$-rectifiable subset of $\R^N$:
 \beq \label{2.3}
 u(\Om) \ = \ \bigcup_{i=1}^\infty \nu^i \big(f(B_i)\big)\ \sub \ \R^N.
 \eeq
 \et
Theorem \ref{th1}  is optimal. Without extra assumptions, there may not exist any $f$ globally defined on $\Om$ and $u(\Om)$ may bifurcate without being given by a single-valued curve $\nu$ for which $u=\nu \circ f$ (Corollary \ref{cor2}, Example \ref{ex1}). Theorem \ref{th1} has been motivated by the rigidity results of Rindler in \cite{R1, R2}. Actually, we extend a part of his result from constant rank-one tensors $\xi \ot w$ to variable rank-one $\xi(x) \ot w(x)$ tensor fields. When compared to the rigidity results known in the literature (see e.g.\ Kirchheim \cite{Ki}), it is somewhat surprising in that most rigidity phenomena appear for rank greater than 2. The idea is as follows: if $Du =\xi \ot w$, then since $\Curl(Du)\equiv 0$, we invoke Poincar\'e's lemma to write $w=Df$ for a scalar $f$ and we also show that $\rk(D \xi)\leq 1$ . Then, we employ geodesic flows, Riemannian exponential maps and a curvilinear extension of ``De Giorgi-type'' arguments to show that $\xi$ and $f$ locally have the same level sets and hence $\xi =\dot{\nu} \circ f$. 

It seems that the natural setting for Theorem \ref{th1} is that of Lipschitz maps. Indeed, we provide such an extension in Theorem \ref{th4}. Yet, this does not follow by a direct approximation argument and substantial complications arise. The problem is that the Rank-One property is \emph{not invariant under mollification}: the mollification may ``fatten'' and its Hausdorff dimension may increase (Remark \ref{rem1}). We remedy this problem by imposing an extra approximation assumption.

In Section \ref{section4} we focus on the general system \eqref{1.7}. We motivate our results by observing that \eqref{1.1} is quasilinear and degenerate elliptic, that is, for
\beq
\A_{\al i \be j}(P)\, := \, P_{\al i} P_{\be j} +|P|^2[P]_{\al \be}^\bot \de_{ij}
\eeq
we can rewrite the $\infty$-Laplacian \eqref{1.1} as $\A(Du)_{\al i \be j} D_{ij}^2u_{\be}=0$ and $\A$ satisfies the symmetry condition and the Legendre-Hadamard condition:
\begin{align} 
 \A_{\al i \be j} \, &=\, \A_{\be j \al i},\\
\A_{\al i \be j}\, \eta_\al a_i \, \eta_\be a_j \, &\geq \, 0, \ \ \ \eta \in \R^N,a\in \R^n.
\end{align}
However, the general system \eqref{1.7} is \emph{not} degenerate elliptic since $[H_P]^\bot \!$ and $H_{PP}$ are symmetric but if $N\geq 2$ their product may \emph{not commute}, not even when $H$ is strictly convex on $\R^{N \by n}$. For $N=1$, though, Aronsson's equation $H_{P_i}H_{P_j}D^2_{ij}u=0$ is trivially degenerate elliptic. In Theorem \ref{th5} we characterise the Hamiltonians which lead to elliptic systems as the ``geometric'' ones which depends on $Du$ via the Riemannian metric $Du^\top \! Du$ on $u(\Om)\sub \R^N$, that is when $H(P) = h \big(\frac{1}{2}P^\top\! P\big)$. In the case of $\De_\infty$, we have $h(p)=\tr (p)$. In dimensions $n\leq 3$, this is a complete equivalence. However, if $n\geq 4$ complicated structures in the higher order tensors $H_{P...P}$ appear and a necessary extra assumption is required for the full equivalence. Without it, $H$ can be written in this form up to an $O(|P|^4)$ correction. In the case $n=1$, we deduce that $H$ is radially symmetric. This is very restrictive, but should be compared with the rigidity of Lipschitz extensions for maps in Kirszbraun's theorem (see e.g.\ Federer \cite{F}, p.\ 201), in contrast to the flexibility of scalar Lipschitz extensions. 

In this paper we also tackle two more independent topics related to the study of solutions to our system \eqref{1.7}. Iin Section \ref{section5} we focus on the 1D case for $n=1\leq N$ and we study the ODE system arising from Hamiltonian $H\in C^2(\R \by \R^N \by \R^N)$ depending on all arguments $H=H\big(x,u(x),u'(x)\big)$. The 1D case of vectorial Calculus of Variations in $L^\infty$ provides an important model for Data Assimilation \cite{K9}. We first formally derive the system in the limit as $p \ri \infty$ of the Euler-Lagrange equations of the respective $L^p$-functional (equation \eqref{5.8}). By imposing the condition of radial dependence in $u'$, we obtain the degenerate elliptic version of the system:
\beq \label{1.21}
A_\infty u  \ = \ |u'|^2\Big(h_p u''\, - \, \bR_{u'}h_\eta\Big) \, + \, h_x u'\ = \  0.
\eeq
Here $h\equiv h\big(\cdot,u,\frac{1}{2}|u'|^2\big)$, the arguments of $h$ are $(x,\eta,p)$ and $\bR_{u'}$ is the \emph{reflection} operator with respect to the normal hyperplane $[u']^\bot$. 
\[
\underset{\text{Figure 4.}}{\includegraphics[scale=0.2]{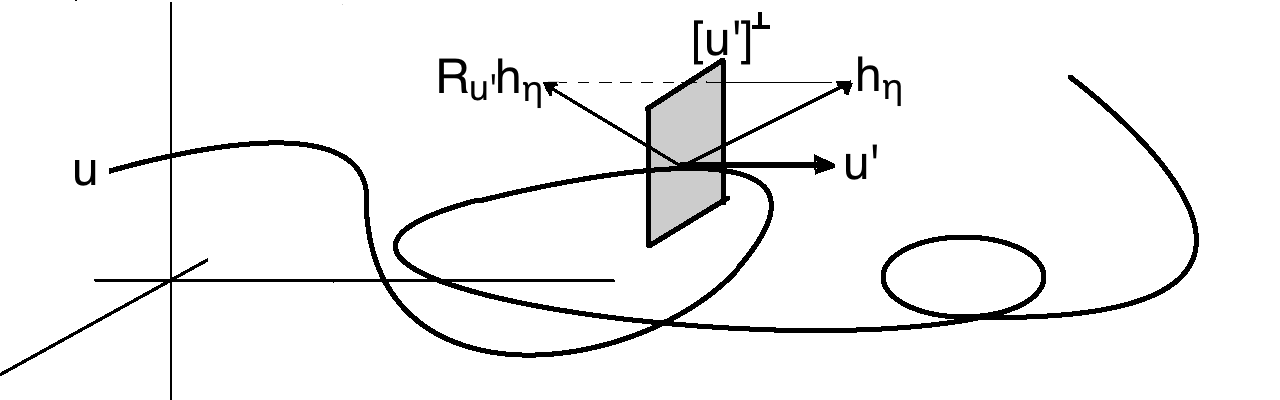}}
\]
We note that although $ \bR_{u'}$ is discontinuous at critical points, in this case the coefficients of \eqref{1.21} are continuous. In Theorem \ref{th6} we study existence, uniqueness and $W^{2,\infty}_{loc}(\R)^N$ regularity of solutions to the initial value problem for \eqref{1.21}.

Finally, motivated by Aronsson's paper \cite{A6}, in Section \ref{section6} we analyse the class of solutions to \eqref{1.1} of the radial form $u=\rho^k f(k \theta)$ for $k>0$ and $f$ a curve in $\R^N$. Interestingly, in Proposition  \ref{Pr1} we prove that such solutions are very rigid since their image is contained in either an affine line or an affine plane.

We conclude this long introduction with some related results known in the literature. In \cite{K4} we identified the \emph{variational principle} characterising $\infty$-Harmonic maps for the model functional $E_\infty(u,\Om)=\|Du\|_{L^\infty(\Om)}$. Surprisingly, the apt notion \emph{is not} the obvious extension of Aronsson's notion in higher dimensions but instead ``\emph{a rank-one absolute minimal coupled by $\infty$-minimal area}". For details see  \cite{K4}. In \cite{K5} we extended the results of \cite{K3}, \cite{K4} to the subelliptic setting. In \cite{K7}, among other things, we proved that the Dirichlet problem for the $\infty$-Laplacian 
\[
\left\{
\begin{array}{l}
\De_\infty u \, =\, 0,\ \ \ \text{ in }\mB^*, \\
\  u(x)\, =\, x,\ \ \ \text{ on }\p \mB^*,
\end{array}
\right.
\]
surprisingly, has \emph{infinitely many} smooth solutions $u : \mB^* \sub \R^n \larrow \R^n$ on the punctured unit ball $\mB^*=\{x:0<|x|<1\}$, for all $n\geq 2$. The crucial observation is that smooth solutions the differential inclusion
\[
Du(\Om) \sub \mathcal{K}\ ,\ \ \ \mathcal{K}\, :=\, \Big\{ P \in \R^{n\by n}  :\, |P|=1,\ \det (P)>0 \Big\}
\]
are $\infty$-Harmonic. In words, smooth solutions of the vectorial Eikonal equation which are local diffeomorphisms solve the $\infty$-Laplacian. For details see \cite{K7}. It is worth mentioning that \emph{the maximum principle we establish herein is not a comparison principle and does not imply uniqueness}. Ou, Troutman and Wilhelm in \cite{OTW} and Wang and Ou in \cite{WO} studied the ``tangential part" of \eqref{1.1}. Sheffield and Smart in \cite{SS} used the nonsmooth operator norm on $\R^{N \by n}$ as their Hamiltonian and derived a very singular variant of \eqref{1.1} which governs the so-called ``tight maps", that is \emph{vectorial optimal Lipschitz extensions}. Our theorem \ref{th3} relates to an analogous phase separation of tight maps observed in \cite{SS}. Capogna and Raich in \cite{CR} used the dilation $K(\P) = \frac {|P|^n}{\textrm{det}(P)}$ as Hamiltonian on $\R^{n \by n}$ and developed an $L^\infty$ variational approach to optimise \emph{Quasiconformal maps}. They derived and studied a special important case of \eqref{1.7}. Their results have been advanced by the author in \cite{K8}.
In the light of our general Theorem \ref{th5}, it is not a coincidence that all Hamiltonians known in the literature depend on the gradient via the Riemannian metric $Du^\top\! Du$.

\subsection{Preliminaries.} \label{Preliminaries} Throughout this paper we reserve $n,N \in \N$ for the dimensions of Euclidean spaces and $\mS^{N-1}$ denotes the unit sphere of $\R^N$. Greek indices $\al, \be, \ga,... $ run from $1$ to $N$ and Latin $i,j,k,...$ form $1$ to $n$. The summation convention will always be employed in repeated indices in a product. Vectors are always viewed as columns. Hence, for $a,b\in \R^n$, $a^\top b$ is their inner product and $ab^\top$ equals $a \ot b$. If $u=u_\al e_\al :  \Om \sub \R^n \larrow \R^N$ is a map, the gradient matrix $Du$ is viewed as $D_i u_\al e_\al \ot e_i : \Om \larrow \R^{N \by n}=\R^N \ot \R^n$ and the Hessian tensor $D^2u$ as $D^2_{ij} u_\al e_\al \ot e_i \ot e_j: \Om \larrow \R^N \ot \mS(\R^n)$. If $V$ is a vector space, then $\mS(V)$ denotes the symmetric linear operators $T : V \larrow V$ for which $T=T^\top$ and $\mS(V)^+$ the subset of nonegative ones. The Euclidean norm on $\R^{N \by n}$ is $|P|=(P_{\al i}P_{\al i})^{\frac{1}{2}} = (\tr (P^\top P))^{\frac{1}{2}}$. If $F\in C^\infty( \R^{N \by n})$ is a function and we denote the standard basis elements of $\R^{N \by n}$ by $e_{\al i}:=e_{\al} \ot e_i$, then its $q$-th order derivative tensor $F_{P...P}$ at $P_0$
\beq
F_{P...P}(P_0)=F_{P_{\al_1 i_1}...P_{\al_q i_q}}(P_0) e_{\al_1 i_1}\ot ...\ot e_{\al_q i_q} 
\eeq 
is viewed as a multilinear map $\ot^{(q)}\R^n \larrow  \ot^{(q)}\R^N$,  or equivalenly as an element of $\ot^{(q)}(\R^{N \by n})$. Here ``$\ot^{(q)}$'' is the $q$-fold tensor product.  Hence, $F_{P...P}$ is a map $\R^{N \by n} \larrow  \ot^{(q)}(\R^{N \by n})$. We will say that a $q$-th order tensor $C \in \ot^{(q)}(\R^{N \by n})$ is fully symmetric in all its arguments when
\beq
C_{...\al i... \be j...}\ = \ C_{...\al j... \be i...} \ = \ C_{... \be j...\al i...}.
\eeq
We also introduce the following \emph{contraction operation} for tensors which extends the inner product $P:Q=\tr(P^\top Q)=P_{\al i}Q_{\al i}$ of $\R^{N \by n}$. For, if $C \in \ot^{(q)}(\R^{N \by n})$ and $A \in \ot^{(p)}(\R^{N \by n})$ with $p\leq q$, we define $C:A \ \in \ \ot^{(q-p)} (\R^{N \by n})$ by
\beq \label{1.24}
(C:A)_{\al_{q} i_{q} ... \al_{p+1} i_{p+1}} \ :=\ C_{\al_q i_q ... \al_1 i_1}  A_{\al_{p} i_{p} ... \al_1 i_1}. 
\eeq
Let now $P : \R^n \larrow \R^N$ be linear map. Upon identifying linear subspaces with orthogonal projections on them, we have the split $\R^N=[P]^\top \oplus [P]^\bot$ where $[P]^\top$ and $[P]^\bot$ denote range of $P$ and nullspace of $P^\top$ respectively. Hence, if $\xi \in \mS^{N-1}$, then $[\xi]^\bot$ or simply $\xi^\bot$ is (the projection on) the normal hyperplane $I-\xi \ot \xi$. Consequently, the $\infty$-Laplacian \eqref{1.1} in index form reads
\beq  \label{1.a}
 D_i u_\al  D_j u_\be D_{ij}^2 u_\be \, +\, |Du|^2 [Du]_{\al \be}^\bot D^2_{ii} u_\be\ = \ 0
\eeq
and the system \eqref{1.7} becomes
\beq \label{1.b}
\Big(H_{P_{\al i}}H_{P_{\be j}}\, +\, H[H_P]_{\al \ga}^\bot H_{P_{\ga i}P_{\be j}}\Big)(Du)D^2_{ij}u_\be
\ = \ 0. 
\eeq
For convenience we use a different scaling in \eqref{1.a} and \eqref{1.b} and we multiply the normal term of \eqref{1.a} by a factor $2$ which is plausible since \eqref{1.b} consists of two systems normal to each other. Finally, $\mH^k$ denotes the $k$-dimensional Hausdorff measure and for measure theoretic notions we use herein we refer to Simon \cite{S}.

\section{Rigidity of Rank-One maps.} \label{section2}

\subsection{The case of smooth Rank-One maps.} \label{subsection2.1} In this subsection we establish our Geometric Analysis rigidity result in the case of $C^2$ maps.

\BPT \ref{th1}. The implication $(ii) \Rightarrow (i)$ is trivial and the whole proof is devoted to establish the reverse implication $(i) \Rightarrow (ii)$. For, suppose there exist $\xi : \Om \sub \R^n \larrow \R^N$ and $w : \Om \sub \R^n \larrow \R^n$ such that $Du = \xi \ot w$. By replacing $\xi$ by $\xi /|\xi|$ on $\{|\xi|>0\}$ and $w$ by $|\xi|w$ on $\{|\xi|>0\}$, we may pass all the zeros of $Du$ to $w$ and assume that $|\xi| \equiv 1$ on
 \beq \label{2.5}
\Om_0 \ := \ \{|Du|>0\}\ = \ \{|w|>0\}.
\eeq
By differentiating $D_k u_\al = \xi_\al w_k$, we have
 \beq \label{2.6}
D^2_{ij}u_\al \ = \ (D_j \xi_\al ) \, w_i \ +\ \xi_\al (D_j w_i),
\eeq
\beq \label{2.7}
D^2_{ji}u_\al \ = \ (D_i\xi_\al ) \, w_j \ +\ \xi_\al (D_i w_j).
\eeq
Since $u \in C^2(\Om)^N$, the curl of $Du$ vanishes and we have
 \beq \label{2.8}
D^2_{ij}u_\al \ = \ D^2_{ji}u_\al.
\eeq
Hence, by \eqref{2.6}, \eqref{2.7}, \eqref{2.8},
 \beq \label{2.9}
\ (D_j \xi) \, w_i \ - \ (D_i \xi) \, w_j \ = \ \xi (D_i w_j \ -  \ D_j w_i).
\eeq
Since $|\xi|^2=1$ on $\Om_0$, we have $D_k \xi ^\top \xi =0$ thereon. Hence, the two sides of \eqref{2.9} are normal to each other. By applying the projections $\xi \ot \xi$ and $[\xi]^\bot = I -\xi \ot \xi$, \eqref{2.9}  decouples on $\Om_0$ to
 \beq \label{2.10}
\Curl (w)_{ij}\ = \ D_i w_j \ -  \ D_j w_i \ \equiv \ 0,
\eeq
\beq \label{2.11}
\ (D_j \xi) \, w_i \ - \ (D_i \xi) \, w_j \ \equiv \ 0.
\eeq
By \eqref{2.10}, the curl of $w : \Om_0 \sub \R^n\larrow \R^n$ vanishes and by \eqref{2.5} $w\equiv 0$ on $\Om \set \Om_0$. Hence, since $\Om$ is contractible, by Poincar\'e 's Lemma $w$ can be represented by the gradient of a scalar function $f \in C^2(\Om)$: $w=Df$. By \eqref{2.11}, for all $i,j \in \{1,...,n\}$ for which $\{w_i \neq 0\} \cap \{w_j \neq 0\} \neq \emptyset$, we have
\beq \label{2.12}
\frac{D_j \xi_\al}{w_j} \ = \ \frac{D_i \xi_\al}{w_i}. 
\eeq
By \eqref{2.12}, the quotient ${D_k \xi_\al}/{w_k}$ is independent of $k$. Hence, we may define 
\beq \label{2.13}
\eta\ := \ \frac{D_k \xi}{w_k} \ \ : \ \ \{w_k \neq 0\} \sub \Om_0 \, \larrow\, \R^N. 
\eeq
By \eqref{2.12}, $\eta$ is well defined on all of $\Om_0$ since $\cup_1^n \{w_k \neq 0\} $ is an open cover of $\Om_0 = \{|w|> 0\}$ and on the overlaps the different expressions coincide. By \eqref{2.13}, we have $D_k \xi_\al = \eta_\al w_k$ on $\{w_k \neq 0\}$. Actually, this extends to the whole of $\Om_0$ since by \eqref{2.11} we get $D_k \xi =0$ whenever $w_k =0$. Thus,
 \beq \label{2.14}
D\xi \ = \ \eta \ot Df \ , \ \ \text{ on }\Om_0 , 
\eeq
and also $\eta$ is normal to $\xi$, since $\eta^\top \xi = \frac{1}{w_k}D_k(\frac{1}{2}|\xi|^2) =0$, on $\{w_k \neq 0\}$. We now employ \eqref{2.14} to show that in a certain local sence $\xi$ and $f$ have the same level sets.

Fix $\al \in \{1,...,N\}$ and set
\begin{align} 
A\ & :=\ \Om_0 \cap \{|\eta_\al |>0\}, \label{2.15}\\
g\ & := \ \xi_\al\ \ , \ \ \la\ := \ \eta_\al. \label{2.17}
\end{align}
We then obtain 
\beq \label{2.18}
 Dg \ = \ \la \, Df\ , \ \ \text{ on }A,
 \eeq
while $|Dg|>0$ and $|\la|>0$ on $A$. \eqref{2.18} says that the level hypersurfaces $\{f=f(x)\}$ and $\{g=g(x)\}$ passing through $x$ have for all $x\in A$ the same tangent spaces:
 \beq \label{2.19}
[Dg]^\bot \ =\ [Df]^\bot \ =\ I -\frac{Df}{|Df|} \ot \frac{Df}{|Df|} .
 \eeq
Consider the level hypersurfaces of $f$, $g$ as Riemannian submanifolds of $A$ with the induced metrics from $\R^n$. Since covariant derivatives coincide with tangential projections of derivatives in $\R^n$, the geodesic equations for $\chi$, $\psi$ with initial conditions $\chi(0)=\psi(0)=x \in A$ and $\dot{\chi}(0)=\dot{\psi}(0)=e \in [Df(x)]^\bot =[Dg(x)]^\bot$ are
\beq
\left\{
\begin{array}{l}
[Df(\chi(t))]^\bot \ddot{\chi}(t)\ = \ 0, \ \ t>0, \ms\\
 \chi(0) = x,\ \ \dot{\chi}(0) =e,
\end{array}
\right.
\eeq
\beq
\left\{
\begin{array}{l}
[Dg(\psi(t))]^\bot \ddot{\psi}(t)\ = \ 0, \ \ t>0, \ms\\
 \psi(0) = x,\ \ \dot{\psi}(0) =e.
\end{array}
\right.
\eeq
Since $[Dg]^\bot \equiv [Df]^\bot $, $\chi$ and $\psi$ satisfiy the same ODEs with the same initial conditions. Hence, by uniqueness, $\chi \equiv \psi$. Consequently, the esponential maps $\exp^f_x$ and $\exp^g_x$ of $\{f=f(x)\}$ and $\{g=g(x)\}$ coincide and hence $(\exp^g_x)^{-1} \circ \exp^f_x$ equals the identity their common geodesically convex neighbourhod centered at $x$. Hence, the level hypersurfaces of $f,\ g$ within $A$ coincide, but perhaps they are at different heights. Cover $A$ by countably many balls whose radii are small enough to guarrantee that the intersections of the level sets of $f$, $g$ with each ball are connected. 
 \[
\underset{\text{Figure 5.}}{\includegraphics[scale=0.25]{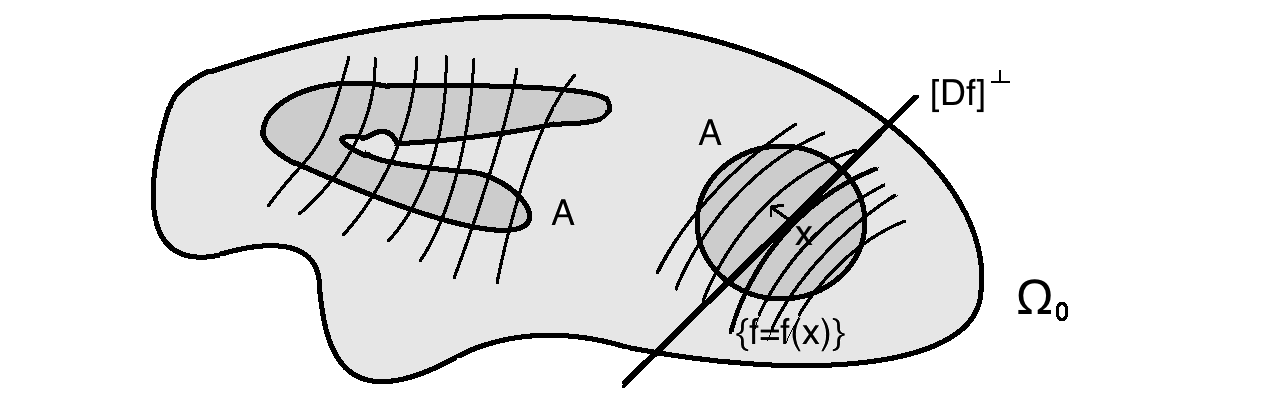}}
\]
Using this cover, we decompose $A$ to a partition of connected Borel sets by writting $A=\cup_1^\infty A_i$, where each $A_i$ equals an open subset of the ball of the cover with possibly some boundary portion. Then, for each $t\in \R$ and each $i \in \N$ there is a unique $\rho^i(t) \in \R$ such that $\{f=t\}$ equals $\{g=\rho^i(t)\}$ locally within $A_i$. Hence, there exists a unique bijection $\rho^i : f(A_i)\sub \R \larrow g(A_i) \sub \R$ such that
 \beq \label{2.22}
\{g=\rho^i(t)\} \ = \ \{f=t\} \ = \ \{\rho^i \circ f=\rho^i(t)\},
\eeq
within $A_i \sub \Om_0$. Equivalently,
 \beq \label{2.23}
 g\ = \ \rho^i \circ f\ , \ \ \text{ on }A_i, \ \ i \in \N.
\eeq
We extend $\rho^i$ from $ f(A_i)$ to $\R$ by zero. 
\[
\underset{\text{Figure 6.}}{\includegraphics[scale=0.23]{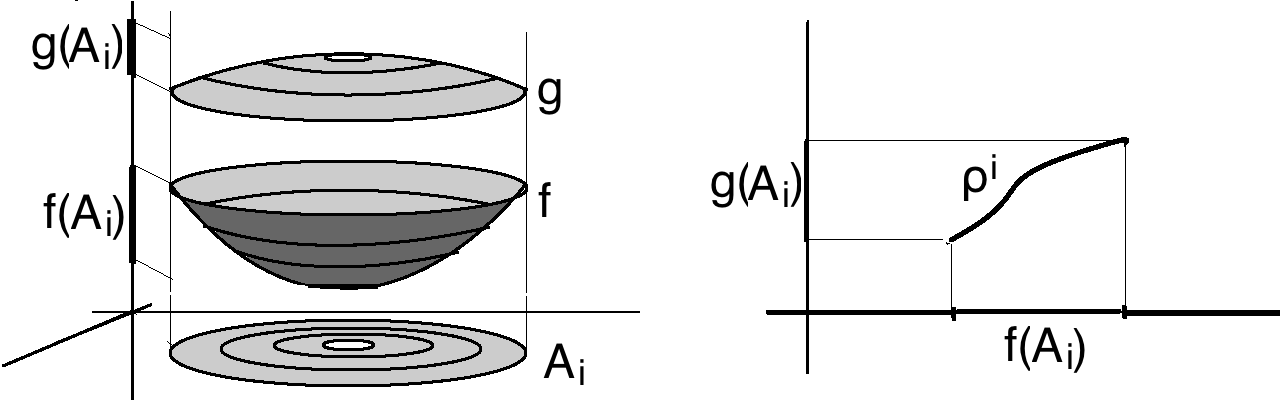}}
\]
On $\Om_0 \set A = \Om_0 \set \cup_1^\infty A_i$, we have $Dg \equiv 0$. Hence, there exists a constant function $\rho^0 : f(\Om_0 \set A) \sub \R \larrow \R$ such that
 \beq \label{2.24}
 g\ = \ \rho^0 \circ f\ , \ \ \text{ on }\Om_0 \set \cup_1^\infty A_i.
\eeq
We extend $\rho^0$ by zero on $\R$ as well. By recalling \eqref{2.15}  and \eqref{2.17}, we have shown that for any $\xi_\al$, $1\leq \al \leq N$, there exists a partition of $\Om_0$ to disjoint connected Borel sets $A_i^\al$ where each $A_i^\al$ equals an open set with possibly some boundary portion and also their complement $A^\al_0 := \Om_0 \set \cup_1^\infty A^\al_i$. There also exist functions $\rho_\al^i : \R \larrow \R$ such that 
\beq \label{2.25}
 \xi_\al\ = \ \rho^i_\al \circ f\ , \ \ \text{ on }A^\al_i, \ \  i=0,1,2,...\ .
\eeq
Hence, by recalling that $|\xi|\equiv 1$ on $\Om_0$, there exists a partitition of $\Om_0$ to connected Borel sets $\{B_i\}_{i \in \N}$ which are intersections of the $A_i$'s and respective bounded curves $\mu^i :\R \larrow \{0\}\cup \mS^{N-1} \sub \R^N$  which satisfy
\beq \label{2.27}
|\mu^i|\equiv 1 \text{ on } f(B_i)\ , \ \ \mu^i \equiv 0 \text{ on }\R \set f(B_i) ,
\eeq
and are such that 
\beq \label{2.26}
 \xi\ = \ \mu^i \circ f\ , \ \ \text{ on }B_i , 
\eeq
for all $i \in \N$. We set
\beq \label{2.28}
\nu^i(t)\ := \ \int_0^t \mu^i(s)\, ds\ , \ \ i \in \N.
\eeq
Then, by \eqref{2.27} we have that $\nu^i \in W^{1,\infty}_{loc}(\R)^N$, while $|\dot{\nu}^i|\equiv 1 \text{ on the interval } f(B_i)$ and also $\dot{\nu}^i \equiv 0 \text{ on }\R \set f(B_i)$. By \eqref{2.26} we have
\beq \label{2.29}
\xi \ = \ \dot{\nu}^i \circ f\ , \ \ \text{ on }B_i .
\eeq
Hence, \eqref{2.29} implies
\begin{align}
Du  \ = \ \xi \ot w \  = \ (\dot{\nu}^i \circ f) \ot Df \ = \ D( {\nu}^i \circ f),
\end{align}
on $B_i$. Thus, $u = {\nu}^i \circ f$ on each $B_i\sub \Om_0$, up to an additive constant. By taking difference quotients in \eqref{2.29}, comparing with \eqref{2.14} and passing to limits, we obtain
\beq
D\xi \ = \ (\ddot{\nu}^i \circ f )\ot Df, 
\eeq
and hence $\ddot{\nu}_\al^i \circ f = D_k\xi_\al\, D_k f$, on $B_i$. Thus, $\ddot{\nu}^i$ exists on $f(B_i)\sub \R$ and is interpreted as $1$-sided at the endpoints of this interval in case it is not open.
Since $Du=0$ and $Df =0$ on $\p(\Om _0) \cap \Om$, we can extend the partition $\cup_1^\infty B_i$ of $\Om_0$ to $\overline{\Om_0} \cap \Om$ and further extend the families $\{B_i\}_{i \in \N}$ and $\{\nu^i\}_{i \in \N}$ by attaching the limit values and setting
\begin{align}
B_0\ & :=  \ \Om \set \overline{\Om_0},\\
\nu^0\ & := \ u\big|_{\Om \set \overline{\Om_0}} \ = \ \text{const}.
\end{align}
Hence, since $u=\nu^i \circ f$ on each $B_i$ of the partition $\cup_0^\infty B_i=\Om$, we conclude that $u$ is $1$-rectifiable and the image $u(\Om)$ equals a union of images of Lipschitz curves:
 \beq
 u(\Om) \ = \ \bigcup_{i=1}^\infty \nu^i \big(f(B_i)\big).
 \eeq
The theorem follows.              \qed

\ms

As we have already mentioned in the Introduction, an extra assumptions is required in order to deduce that a rank-one map $u$ has the form $u=\nu \circ f$ for a unique single-valued unit speed curve $\nu$. This assumption guarrantees ``low complexity'' for the direction field $\xi$.

\bcor[Strong Rigidity of Rank-One maps] \label{cor2}
Suppose $\Om \sub \R^n$ is open and contractible and $u : \Om \sub \R^n \larrow \R^N$ is in $C^2(\Om)^N$. Consider the following statements:

\ms 

\noi (i) $u$ is a strictly Rank-One map, that is $\rk(Du) = 1$ on $\Om$ or equivalently there exist $C^1$ maps $\xi : \Om \larrow \R^N \set \{0\}$ and $w : \Om \larrow \R^n \set \{0\}$ such that $Du = \xi \ot w$. Moreover, the following condition holds
 \beq \label{2.4}
 E\ := \ \Om \cap \left(\bigcup_{\al =1}^N \p \big\{ |D\xi_\al|>0\big\}\right) \ = \ \emptyset.
\eeq

\noi (ii) 
 $u$ equals the composition of a single curve $\nu \in W^{1,\infty}_{loc}(\R)^N$ with a scalar function $f \in C^2(\Om)$, without critical points that is $u=\nu \circ f$ with $|\dot{\nu}| \equiv 1$ on $f(\Om)$, $\dot{\nu} \equiv 0$ on $\R \set f(\Om)$. Moreover, $Du \ =\ (\dot{\nu} \circ f ) \ot Df$ on $\Om$ and $u(\Om)$ is $1$-rectifiable, equal to $\nu(f(\Om))$.

\ms

Then, (i) implies (ii) and also (ii) implies that $u$ is a strictly rank-one map, that is assertion (i) without \eqref{2.4}.
\ecor

\BPCOR \ref{cor2}. In the setting of the proof of Theorem \ref{th3}, if in addition the set $ E$ given by \eqref{2.4} is empty and moreover $\rk(Du)>0$ on $\Om$, then for all $\al \in \{1,...,N\}$, either $D\xi_\al$ does not vanish anywhere inside $\Om_0=\Om$ or it is identically constant. In both cases, the previous set $A$ is connected and coincides with $\Om$. Hence, the curve $\nu$ constructed is unique and consequently $u=\nu \circ f$ with $|\dot{\nu}| \equiv 1$ on $f(\Om)$ and $\dot{\nu} \equiv 0$ on $\R \set f(\Om)$. The reverse implication is obvious.
\qed

\ms

\bex \label{ex1} The additional assumption \eqref{2.4} of Corollary \ref{cor2} is necessary in order to obtain $u=\nu \circ f$. It reduces the complexity of $\xi$ and leads to the avoidance of bifurcations in the curve $\nu$. For, let $u : \R^2 \larrow \R^2$ be given by
\[
u(x)\ := \ 
\left\{
\begin{array}{l} \big(+f^4(x),f(x)\big)^\top,\ \ \text{ on }\{f>0\}\cap\{x_1>0\},\ms\\
\big(-f^4(x),f(x)\big)^\top,\ \ \text{ on }\{f>0\}\cap\{x_1<0\},\ms\\
(0,f(x))^\top, \hspace{35pt} \text{ on }\{f\leq 0\},
\end{array}
\right.
\]
where
\[
f(x)\ :=  \ 1 \, - \, |x-e_1|^2|x+e_1|^2.
\]
Then, $u$ can not be written as $u=\nu \circ f$ for a single-valued curve $\nu$ since  the unique  $\nu$ bifurcates and has two branches: $\nu^\pm (t)=\big(\pm t^4\chi_{(0,\infty)}(t),t\big)^\top $.
 \[
\underset{\text{Figure 7.}}{\includegraphics[scale=0.23]{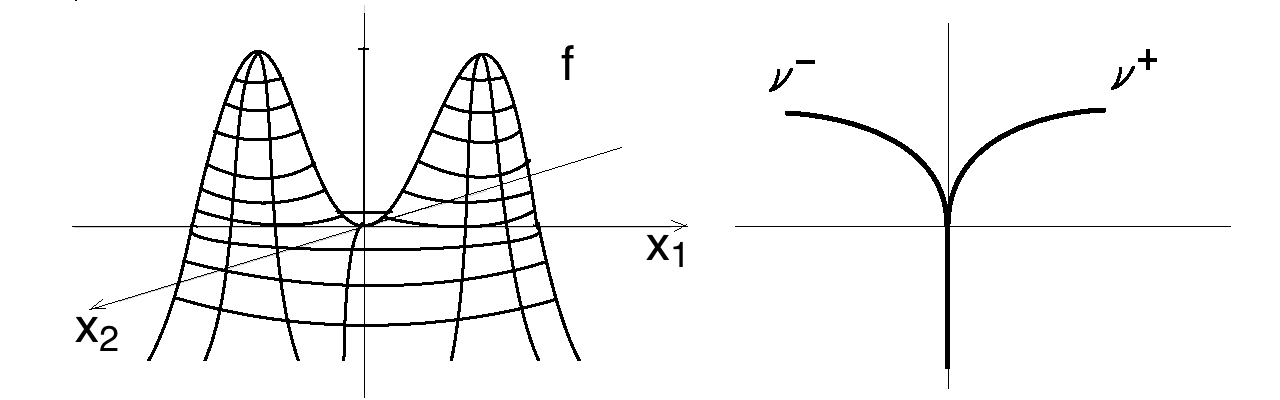}}
\]
\eex

\subsection{Extension to Lipschitz Rank-One maps.} In this subsection we extend Theorem \ref{th1} to the Lipschitz setting. As we have already explained, this does not follow by a standard mollification argument and an additional approximation property is required, which we introduce as assumption.

 \bt[Rigidity of Lipschitz Rank-One maps] \label{th4} Suppose $\Om \sub \R^n$ is open, bounded and contractible and $u : \Om \sub \R^n \larrow \R^N$ is in $W^{1,\infty}(\Om)^N$. 

We moreover assume that there exists a family $\{V^\e \}_{\e>0}$ of rank-one smooth tensor fields in $C^\infty (\Om)^{Nn}$ where each $V^\e$ is curl-free (that is $\rk(V^\e)\leq 1$ and also $D_j V^\e_{ \al i} - D_iV^\e_{ \al j} =0$) such that 
\beq \label{2.35}
V^\e \, \weakstar \, Du \text{ \ in } L^\infty(\Om)^{Nn} \text{ and\ \ $V^\e  \larrow Du$ a.e.\ on } \Om,\ \  \text{ as } \ \e \ri 0 .
\eeq
Then, the following are equivalent:

\ms 

\noi (i) $u$ is a Rank-One map, that is $\rk(Du)\leq 1$ a.e.\ on $\Om$ or equivalently there exist $L^\infty$ vector fields $\xi : \Om \larrow \R^N$ and $w : \Om \larrow \R^n$ such that $Du = \xi \ot w$ a.e.\ on $\Om$.

\ms 

\noi (ii) There exists $f \in W^{1,\infty}(\Om)$, a partition $\{B_i\}_{i \in \N}$ of $\Om$ to measurable sets which covers it a.e., that is $\big|\Om\set (\cup_1^\infty B_i)\big|=0$ and Lipschitz curves $\{\nu^i \}_{i \in \N} \sub W^{1,\infty}_{loc}(\R)^N$ such that on each $B_i$ $u$ equals the composition of $\nu^i$ with $f$:
 \beq \label{2.36}
 u\ = \ \nu^i \circ f\ , \ \ \text{ on }B_i \sub \Om.
  \eeq 
Moreover, $\|\dot{\nu}^i\|_{L^\infty(\R)} \leq 1$ and $\dot{\nu}^i = 0$ a.e.\ on $\R \set f(B_i)$. Also, 
 \beq \label{2.37}
  Du \ =\ (\dot{\nu}^i \circ f ) \ot Df \ , \ \ \text{ a.e. on }B_i \sub \Om,
 \eeq
and the image $u(\Om)$ is an $1$-rectifiable subset of $\R^N$:
 \beq \label{2.38}
\mH^1\left( u(\Om) \ \set \ \bigcup_{i=1}^\infty \nu^i \big(f(B_i)\big)\right)\ = \ 0.
 \eeq
 \et

\begin{remark} \label{rem1}
The extra approximation assumption \eqref{2.35} of Theorem \ref{th4} requires that $Du$ is in the intersection of the weak$^*$ and the pointwise closures in $ L^\infty(\Om)^{Nn}$ of the cone which consists of smooth rank-one curl-free tensor fields. Such an assumption is superfuous if either $\xi$ or $w$ is identically constant, since mollification of $Du =\xi \ot w$ produces the desired approximations $V^\e$. 

Generally, however, all standard mollification methods average at each point contributions from nearby points. As a result, if such a ``partial affinity'' of $u$ fails to hold and both $\xi$ \emph{and} $w$ vary, the range $u(\Om)$ may ``fatten'' and the mollification of $u$ may not be rank-one any more. Unfortunately, we have not been able neither to verify the necessity of the assumption nor to construct a proper mollification scheme allowing to drop it. Notwithstanding, this $W^{1,\infty}$-extension is not required for the phase separation theorem of the $\infty$-Laplacian.
\end{remark}

\BPT \ref{th4}. Is suffices to demonstrate the implication $(i)\Rightarrow (ii)$. Suppose $Du=\xi \ot w$ a.e.\ on $\Om$. By a rescaling of the form $Du=(\frac{1}{|\xi|}\xi ) \ot (|\xi|w)$ on $\{|\xi|>0\}$, we may assume that $\xi : \Om_0 \larrow \mS^{N-1}$, where $\Om_0 := \{|Du|>0\} \sub \Om$ and also that $\xi =0$ a.e.\ on $\Om \set \Om_0$. By assumption, we have $\rk(V^\e)\leq 1$ and hence there exist $\xi^\e : \Om\sub \R^n \larrow \R^N$ and  $w^\e : \Om \sub \R^n \larrow \R^n$ such that $V^\e = \xi^\e \ot w^\e$. By an appropriate rescaling inside the products $(\frac{1}{|\xi^\e|} \xi^\e) \ot ({|\xi^\e|}w^\e)$ on $\{|\xi^\e| >0\}$, we may assume that $\xi^\e : \Om_\e \larrow \mS^{N-1}$ where $\Om_\e := \{|V^\e |>0\} \sub \Om$ and also that $\xi^\e \equiv0$ on $\Om \set \Om_\e$. 

We now claim that $\xi^\e \larrow \xi $ and also that $ w^\e \larrow w$ as $\e \ri 0$, both weakly$^*$ in $L^\infty(\Om)$ and also a.e.\ on $\Om$; indeed, there exists $\eta$ such that $\xi^\e \overset{ \phantom{a}_* \ }{\lharpoonup} \eta$ and hence by the $L^1(\Om)^{Nn}$ strong convergence of $\xi^\e \ot w^\e$ which follows by the Dominated Convergence theorem, we have
\beq
w^\e =\, (\xi^\e)^\top (\xi^\e \ot w^\e) \ \weakstar \ \eta^\top (\xi \ot w)\, =\,  (\eta^\top \xi)w,
\eeq
as $\e \ri 0$. Thus, by uniqueness of limits of $\xi^\e \ot w^\e$ we have $[(\eta \ot \eta)\xi]\ot w = \xi \ot w$ a.e.\ on $\Om$ and hence $\xi =\eta$. Since $\Om$ is contractible, by Poincar\'e's lemma, for any $\e>0$ there exists a smooth map $u^\e : \Om \sub \R^n \larrow \R^N$ such that $V^\e$ can be represented as the gradient of $u^\e$: $Du^\e = \xi^\e \ot w^\e$. Moreover, each $u^\e$ is a smooth rank-one map: by Theorem \ref{th1}, there exist scalar functions $f^\e \in C^\infty(\Om)$, partitions of $\Om$ to Borel sets $\{B^\e_i\}_{i \in \N}$ with $\Om=\cup_1^\infty B^\e_i$, families of Lipschitz curves $\{\nu^{i\e}\}_{i \in \N} \sub W^{1,\infty}_{loc}(\R)^N$ with $\|\dot{\nu}^{i\e}\|_{L^\infty(\R)}\leq 1$ and $\dot{\nu}^{i\e}\equiv 0$ on $\R \set f^\e(B^\e_i)$ such that $u^\e = \nu^{i\e} \circ f^\e$ on each $B_i^\e \sub \Om$, while the images $u^\e(\Om)$ are $1$-rectifiable, equal to $\cup_1^\infty \nu^{i\e} (f^\e(B_i^\e))$. 

 We will now show that appropriate normalised shifts of the maps $u^\e$ approximate $u$. Fix a point $\overline{x} \in \Om$ and set $d:= \diam(\Om)$. Since $Du^\e \weakstar Du$ in $L^\infty(\Om)^{Nn}$ as $\e \ri 0$, for all $x,y \in \Om$ and $\e >0$ small we have
\begin{align} \label{2.39}
|u^\e(x) - u^\e(y)| \ \leq \ \big(\|Du\|_{L^\infty(\Om)} + 1\big) |x-y|.
\end{align}
We further normalise $u^\e$ by considering appropriate shifts, denoted again by $u^\e$, such that $u^\e(\overline{x} )= u(\overline{x} )$. By \eqref{2.39}, we have
\begin{align} \label{2.39a}
\|u^\e\|_{L^\infty(\Om)} \ \leq \ d\big(\|Du\|_{L^\infty(\Om)} + 1\big) \ + \ |u(\overline{x})|.
\end{align}
Hence, there exists $v$ such that $u^\e \weakstar v$ in $W^{1,\infty}(\Om)^{N}$ as $\e \ri 0$. We will now show that $u\equiv v$. Since $Du^\e \larrow Du$ a.e.\ on $\Om$, for $\mH^{n-1}$-a.e.\ direction $e\in\mS^{n-1}$,  we have that $Du^\e \larrow Du$ $\mH^1$-a.e.\ on the set $(\overline{x} + \spn [e])\cap \Om =: I$. We fix such an $e$. By Egoroff's theorem, for any $\si \in (0,1)$, there is an $\mH^1$-measurable set $E_\si \sub I$ with $\mH^1(E_\si)\leq \si$ such that $Du^\e \larrow Du$ uniformly on $I \set E_\si$ as $\e \ri 0$. Since $u^\e(\overline{x} )= u(\overline{x} )$, by the $1$-dimensional Poincar\'e inequality, for $\e>0$ small we have
\begin{align} \label{2.40}
\int_{0}^{d} \big|u^\e(\overline{x}+te)-  u(\overline{x}+te) & \big|\, dt\  \leq \ d \int_{0}^{d} \big|Du^\e(\overline{x}+te)e-Du(\overline{x}+te)e\big|\, dt \nonumber\\
& \leq \ d^2 \sup_{I \set E_\si}\big|Du^\e - Du\big| \, + \, d \big(2\|Du\|_{L^\infty(\Om)} + 1\big) \mH^1(E_\si). 
\end{align}
Since $u^\e \larrow v$ in $C^0\big(\overline{\Om}\big)^N$ and $Du^\e \larrow Du$ in $C^0(I \set E_\si)^{Nn}$ as $\e \ri 0$, by passing to the limit in \eqref{2.40} we obtain
\beq 
\label{2.41}
\int_0^d \big|v(\overline{x}+te)-u(\overline{x}+te)\big|\, dt\  \leq \ d\, \big(2\|Du\|_{L^\infty(\Om_R)} + 1\big)\si . 
\eeq
By letting $\si \ri 0$, by \eqref{2.41} we get $u\equiv v$ on $I\sub \Om$. Since this holds for $\mH^{n-1}$-a.e.\ direction $e \in \mS^{N-1}$, we get $u\equiv v$ on $\Om$. Hence, $u^\e \weakstar u$ in $W^{1,\infty}(\Om)^{N}$ as $\e \ri 0$. Since $Df^\e \weakstar w$ in $L^\infty(\Om)^n$ , for $\e>0$ small we have
 \begin{align} \label{2.42}
|f^\e(x) - f^\e(y)| \ \leq \ \big(\|w\|_{L^\infty(\Om)} + 1\big) |x-y|.
\end{align}
We further normalise the family $f^\e$ by considering appropriate shifts denoted again by $f^\e$ such that $f^\e(\overline{x})=f(\overline{x})$. By replacing also each $\nu^{\e i}$ with the translate $\nu^{\e i}\big(\cdot -(f(\overline{x})-f^\e(\overline{x}))\big)$, we do not affect the previous normalisation  $u^\e(\overline{x})=u(\overline{x})$. Consequenly, \eqref{2.42} implies
 \begin{align} \label{2.42a}
\|f^\e\|_{L^\infty(\Om)} \ \leq \ d\big(\|w\|_{L^\infty(\Om)} + 1\big)\ +  \ |f(\overline{x})| .
\end{align}
As a result, there exists an $f$ such that $f^\e \weakstar f$ in $W^{1,\infty}(\Om)$ as $\e \ri 0$.  

Since $\dot{\nu}^{\e i} \circ f^\e=\xi^\e$ on $B^\e_i$ and $\dot{\nu}^{\e i} \circ f^\e = 0$ on $\Om \set B^\e_i$, for $\e,\de >0$ small we have
 \begin{align} \label{2.43}
\big|B^\e_i \triangle B^\de_i\big| \ & = \ \int_{\Om}\big|\chi_{B^\e_i} - \chi_{B^\de_i}\big|\nonumber\\
& = \ \int_{\Om}\big||\dot{\nu}^{\e i} \circ f^\e | - |\dot{\nu}^{\de i} \circ f^\de|\big|\\
& \leq \ \int_{\Om}\big|\xi^\e - \xi^\de\big|. \nonumber
\end{align}
Since $\xi^\e \larrow \xi $ in $L^1(\Om)^N$, for each $i\in \N$ the family $\{B^\e_i\}_{\e>0}$ is Cauchy in measure and hence has a measurable limit $B_i \sub \Om$. Since for all $\e>0$ we have $\Om = \cup_1^\infty B^\e_i$ and $ B^\e_i \cap  B^\e_j =\emptyset$ for $i\neq j$, the limit family forms a cover of $\Om$ except perhaps for a nullset: $\big|\Om \set (\cup_1^\infty B_i)\big|=0$. We recall that we have $u^\e = \nu^{\e i} \circ f^\e$ on $B^\e_i$ and also $\|\dot{\nu}^{\e i}\|_{L^\infty(\R)}\leq 1$ and $\dot{\nu}^{\e i} \equiv 0$ on $\R \set f^\e(B^\e_i)$. 
 \[
\underset{\text{Figure 8.}}{\includegraphics[scale=0.22]{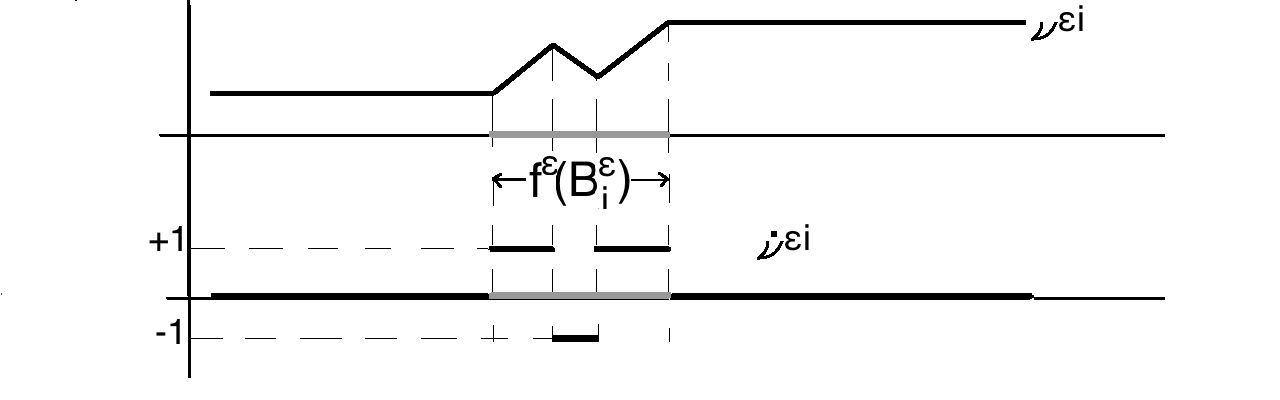}}
\]
Hence, if $\overline{x} \in B^\e_i$, for any $t\in \R$ we have
 \begin{align} \label{2.43a}
|\nu^{\e i}(t)|\ &\leq \ \|\dot{\nu}^{\e i}\|_{L^\infty(\R)} |t-f^\e(\overline{x})|\ + \ |\nu^{\e i}(f^\e(\overline{x}))| \nonumber \\
& = \ |t-f(\overline{x})|\ + \ |u(\overline{x})|. 
\end{align}
If $\overline{x} \not\in B^\e_i$, then $f(\overline{x})$ is in the complement of the interval $f^\e(B^\e_i)$ and since $|\nu^{\e i}|$ is constant on $\R \set f^\e(B^\e_i)$, for any $t\in \R$ we have
 \begin{align} \label{2.43b}
|\nu^{\e i}(t)|\ &\leq \ \|\dot{\nu}^{\e i}\|_{L^\infty(\R)} |t-f(\overline{x})|\ + \ |\nu^{\e i}(f(\overline{x}))|\nonumber \\
 &\leq \ |t-f(\overline{x})|\ + \ \max \big|\nu^{\e i}(\p (f^\e(B^\e_i))\big| \\
& \leq \ |t-f(\overline{x})|\ + \  \|u^\e\|_{L^\infty(\Om)}. \nonumber
\end{align}
As a result, since the family $u^\e$ is uniformly bounded on $\Om$, for each $i \in \N$ the family $\{\nu^{\e i}\}_{\e>0} \sub W^{1,\infty}_{loc}(\R)^N$ has a weak$^*$ limit $\nu^i \in W^{1,\infty}_{loc}(\R)^N$ which satisfies $\|\dot{\nu}^i \|_{L^\infty(\R)}\leq 1$. By passing to the limit as $\e \ri 0$ we get  $u = \nu^{i} \circ f$ on $B_i \sub \Om$ and $\nu^i = 0$ on $\R \set f(B_i)$. Finally, the image $u(\Om)$ is $1$-rectifiable in $\R^N$ and up to an $\mH^1$-nullset of $\R^N$, we have $u(\Om)=\cup_1^\infty \nu^i(f(B_i))$. The theorem follows.                           \qed

\section{The structure of $2$-dimensional $\infty$-Harmonic maps.} \label{section3}

In this section we use the Rigidity Theorem \ref{th1} proved in Section \ref{section2} to analyse the phase separation of classical solutions to \eqref{1.1} when $n=2$ and $N\geq 2$. 

\BPT \ref{th3}. We begin by setting
\begin{align}
\Om_1\ &  := \ \inter \big\{ \rk (Du) \leq 1\big\}, \label{3.5}\\
\Om_2\ & := \ \big\{ \rk (Du) =2 \big\},  \label{3.6}
\end{align}
and let also $\S:=\Om \set(\Om_1 \cup \Om_2)$. Our PDE system \eqref{1.1} decouples to
\begin{align}
Du\, D\Big(\frac{1}{2}|Du|^2\Big)\ = \ 0, \label{3.7}\\
|Du|^2[Du]^\bot \De u \  =\ 0.  \label{3.8}
\end{align}
On $\Om_2$, we have $\rk(Du)=2$ and hence $u\big|_{\Om_2} : \Om_2 \larrow \R^N$ is an immersion. Thus, $Du(x)$ possesses a left inverse $(Du(x))^{-1}$  for all $x\in \Om_2$. Hence, \eqref{3.7} implies
\beq
(Du)^{-1} Du\, D\Big(\frac{1}{2}|Du|^2\Big)\ = \ 0
\eeq
and hence $ D\Big(\frac{1}{2}|Du|^2\Big) =  0$ on $\Om_2$, or equivalently
\beq \label{3.10}
|Du|^2\ = \ \text{const.},
\eeq
on each connected component of $\Om_2$. Moreover, \eqref{3.10} holds on $\S$ as well, the common boundary of $\Om_2$ and $\Om_1$. 
\[
\underset{\text{Figure 9.}}{\includegraphics[scale=0.2]{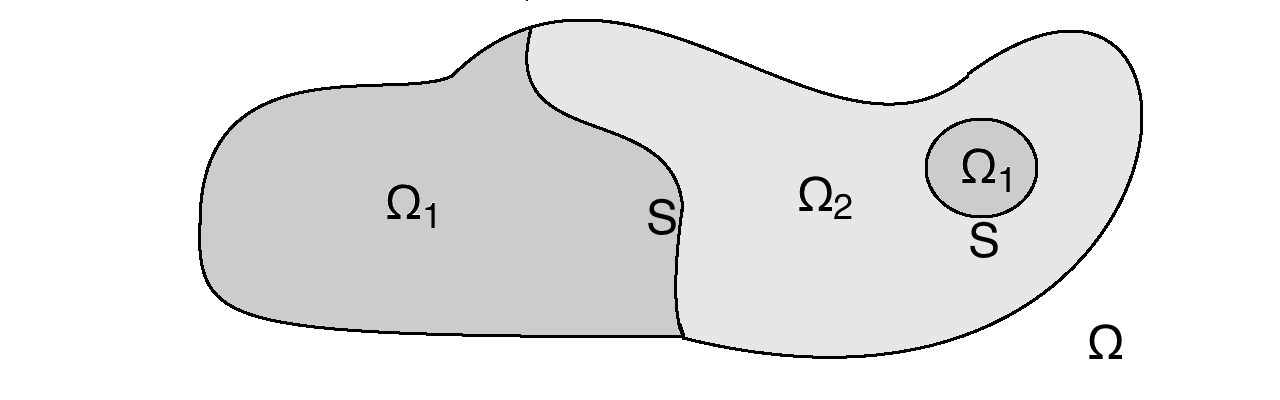}}
\]
On the other hand, on $\Om_1$ we have $\rk(Du)\leq 1$. Hence, there exist vector fields $\xi : \Om_1 \sub \R^2 \larrow \R^N$ and $w : \Om_1 \sub \R^2 \larrow \R^n$ such that $Du =\xi \ot w$. Suppose first that $\Om_1$ is contractible. Then, by the Rigidity Theorem \ref{th1}, there exists a function $f\in C^2(\Om_1)$, a partition of $\Om_1$ to Borel sets $\{B_i\}_{i \in \N}$ and Lipschitz curves $\{\nu^i\}_{i \in \N} \sub W^{1,\infty}_{loc}(\R)^N$ with $|\dot{\nu}^i|\equiv 1$ on $f(B_i)$, $|\dot{\nu}^i|\equiv 0$ on $\R \set f(B_i)$ twice differentiable on $f(B_i)$, such that $u = \nu^i \circ f$ on each $B_i$ and hence $Du = (\dot{\nu}^i \circ f) \ot Df$ on $B_i$. By \eqref{3.7}, we obtain
\begin{align}
\big((\dot{\nu}^i \circ f) \ot Df & \big) \ot  \big((\dot{\nu}^i \circ f) \ot Df\big) : \nonumber \\
 &: \Big[(\ddot{\nu}^i \circ f) \ot Df \ot Df \, +\, (\dot{\nu}^i \circ f) \ot D^2f\Big]\ =\ 0,
\end{align}
on $B_i \sub \Om_1$. Since $|\dot{\nu}^i |\equiv 1$ on $f(B_i)$, we have that $\ddot{\nu}^i$ is normal to $\dot{\nu}^i$ and hence
\beq
\big((\dot{\nu}^i \circ f) \ot Df \big) \ot  \big((\dot{\nu}^i \circ f) \ot Df\big) : \big((\dot{\nu}^i \circ f) \ot D^2f\big)\ =\ 0,
\eeq
on $B_i \sub \Om_1$.  Hence, by using again that $|\dot{\nu}^i |^2\equiv 1$ on $f(B_i)$ we get
\beq
\big(Df \ot Df : D^2f \big) (\dot{\nu}^i \circ f)  \ =\ 0,
\eeq
on $B_i \sub \Om_1$. Thus, $\De_\infty f = 0$ on $B_i$. By \eqref{3.8} and again since  $|\dot{\nu}^i |^2\equiv 1$ on $f(B_i)$, we have $[Du]^\bot =  [\dot{\nu}^i \circ f]^\bot$ and hence
\beq
|Df|^2 \, [\dot{\nu}^i \circ f]^\bot \Div\big((\dot{\nu}^i \circ f) \ot Df\big) \ =\ 0,
\eeq
 on $B_i \sub \Om_1$.  Hence,
\beq
|Df|^2 \, [\dot{\nu}^i \circ f]^\bot \Big((\ddot{\nu}^i \circ f)|Df|^2 \, +\, (\dot{\nu}^i \circ f) \De f\Big) \ =\ 0,
\eeq
on $B_i$, which by using once again  $|\dot{\nu}^i |^2\equiv 1$ gives
\beq \label{3.16}
|Df|^4  (\ddot{\nu}^i \circ f)\ =\ 0,
\eeq
on $B_i$. Since $\De_\infty f = 0$ on $B_i$ and $\Om_1 = \cup_1^\infty B_i$, $f$ is $\infty$-Harmonic on $\Om_1$. Thus, by Aronsson's theorem in \cite{A4}, either $|Df|>0$ or $|Df|\equiv 0$ on $\Om_1$.

If the first alternative holds, then by \eqref{3.16} we have $\ddot{\nu}^i \equiv 0$ on $f(B_i)$ for all $i$ and hence $\nu^i$ is affine on $f(B_i)$, that is $\nu^i(t)=t\xi^i +a^i$ for some $|\xi^i|=1$, $a^i \in \R^N$. Thus, since $u=\nu^i \circ f$ and $u \in C^2(\Om_1)^N$, all $\xi^i$ and all $a^i$ coincide and consequently $u=\xi f + a$, $\xi \in \mS^{N-1}$, where $a \in \R^N$ and $f \in C^2(\Om_1)$.

If the second alternative holds, then $f$ is constant on $\Om_1$ and hence by the representation $u=\nu^i \circ f$, $u$ is piecewise constant on each $B_i$. Since $u \in C^2(\Om_1)^N$ and $\Om_1=\cup_1^\infty B_i$, necessarily $u$ is constant on $\Om_1$. But then $|Du|_{\Om_2}|=|Df|_{\S}|=0$ and necessarily $\Om_2 =\emptyset$. Hence, $|Du|\equiv 0$ on $\Om$, that is $u$ is affine on each of the connected components of $\Om$.

If $\Om_1$ is not contractible, cover it with balls $\{\mB_m\}_{m\in \N}$ and apply the previous argument. Hence, on each $\mB_m$,  we have  $u=\xi^m f^m + a^m$, $\xi^m \in \mS^{N-1}$,  $a^m \in \R^N$ and $f^m \in C^2(\mB_m)$ with $\De_\infty f^m = 0$ on $\mB_m$ and hence either $|Df^m|>0$ or $|Df^m|\equiv 0$. Since $u \in C^2(\Om_1)^N$, on the overlaps of the balls the different expressions of $u$ must coincide and hence  we obtain $u=\xi f + a$ for $\xi \in \mS^{N-1}$,  $a \in \R^N$ and $f \in C^2(\Om_1)$ where $\xi$ and $a$ may vary on different connected components of $\Om_1$. The theorem follows.  \qed

\ms

Theorem \ref{th3} implies a vectorial version of the Maximum Principle when $n=N=2$, which we now prove.

\BPCOR \ref{Cor2}. We begin by observing that \eqref{3.17} is an elegant restatement of the Maximum Principle for all projections $\eta^\top u$ of $u$, that is, when for all $\Om' \Subset \Om$ and all directions $\eta \in \mS^{N-1}$ we have
\beq \label{3.18}
\sup_{\Om'}\eta^\top u \ \leq \ \max_{\p \Om'} \eta^\top u.
\eeq
Indeed, \eqref{3.18} says that $u(\Om')$ is contained in the intersection of all halfspaces containing $u(\p \Om')$. To see \eqref{3.18}, fix $\Om' $ and  $\eta \in \mS^{N-1}$ and let $\Om_1$, $\Om_2$, $\S$ respectively be the constant rank domains and the interface of $u$, as in Theorem \ref{th3}. Suppose that $u=\xi f +a$ on $\Om_1\cup \S$, where $\xi \in \mS^{N-1}$, $a\in \R^N$ and $f\in C^2(\Om_1 \cup \S)$. Then, 
\begin{align} \label{3.19}
|D(\eta^\top u)|\ & \ = \ |\eta^\top Du|\chi_{\Om_2} \ +\  |\eta^\top Du|\chi_{\S \cup \Om_1} \nonumber\\
                        & \ = \ |\eta^\top Du|\chi_{\Om_2} \ +\  |\eta^\top\xi|\, | Df|\chi_{\S \cup \Om_1}.
\end{align}
If $|Df|\equiv 0$ on $\Om_1$, then $\Om_2=\emptyset$ and $u$ is affine. Hence, \eqref{3.18} follows. Suppose now $|Df|>0$ on $\Om_1$. Since $u|_{\Om_2}$ is a local diffeomorphism, we have $|\eta^\top Du|>0$ for all $\eta \in \mS^{N-1}$. 
 \[
\underset{\text{Figure 10.}}{\includegraphics[scale=0.2]{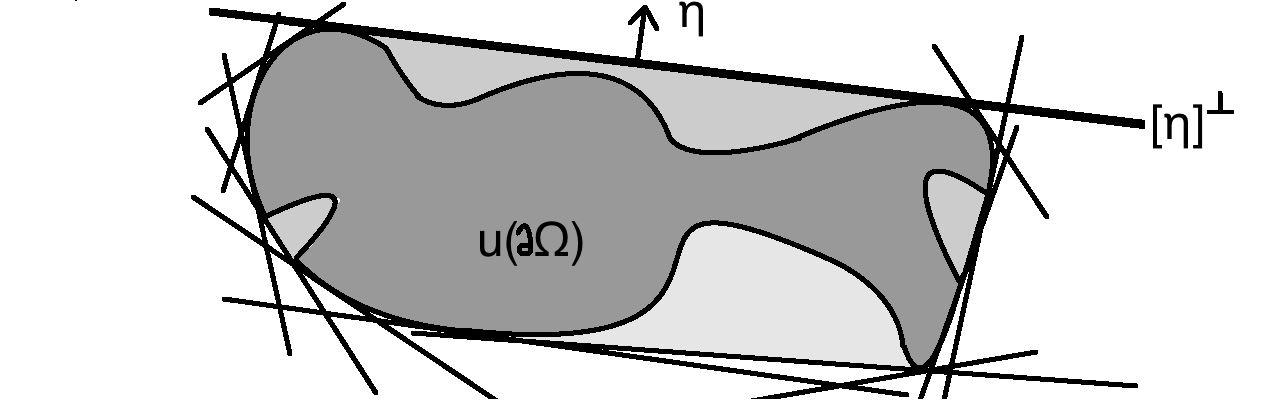}}
\]
Consequently, for all $\eta \in \mS^{N-1} \set [\xi]^\bot$, in view of \eqref{3.19} we have $|D(\eta^\top u)|>0$  on $\Om$. Hence, $\eta^\top u$ has no interior critical points inside $\Om$ and consequently we have
\beq \label{3.20}
\max_{\overline{\Om'}}\eta^\top u \ = \ \max_{\p \Om'} \eta^\top u,
\eeq
for all directions $\eta \not\perp \xi$. By letting $\dist(\eta,[\xi]^\bot) \ri 0$, \eqref{3.20} implies \eqref{3.17}.
\qed

\ms

\section{Characterisation of the class of elliptic PDE systems.} \label{section4}

In this section we focus on the general Aronsson system \eqref{1.7}. As already explained in the introduction, when $N\geq 2$ the normal coefficient $H[H_P]^\bot H_{PP}$ is not symmetric and as a result the system generally is not degenerate elliptic, not even for strictly convex Hamiltonians. In Theorem \ref{th5} below we establish that all ``geometric'' Hamiltonians which depend on $Du$ via the induced Riemannian metric $Du^\top\! Du$ lead to elliptic systems. Moreover, in low dimensions $n\leq 3$ the converse is true as well for (normalised) analytic Hamiltonians with fully symmetric Hessian tensor. When $n\geq 4$, there appear complicated structures in the minors of forth and higher order derivatives and an additional assumption is required. The constructive method of proof reveals that it is necessary. The main idea in the reverse direction is to impose the commutativity relation $[H_P]^\bot H_{PP}=H_{PP}[H_P]^\bot $ and use power-series expansions of $H$ and induction, by a term-after-term blow-up argument along inverse images under $H_P$ of rank-one directions.

\bt[Classification of Hamiltonians leading to elliptic systems \eqref{1.7}] \label{th5} Suppose that $H\in C^2(\R^{N \by n})$ is a non-negative Hamiltonian with $n\geq 1$, $N\geq 2$. Suppose also that  $[H_P(P)]^\bot=[P]^\bot$ on $\R^{N\by n}$. Consider the following statements:

\ms

\noi (i) There exists $h \in C^2\left(\mS(\R^n)^+\right)$, with symmetric gradient $h_p$, such that
 \beq \label{4.1}
H(\P)\ = \ h\Big({\frac{1}{2} }P^\top P\Big).
\eeq

\noi (ii) The system 
\beq \label{4.2}
\A_\infty u \, :=\, \Big(H_P \ot H_P+ H\, [H_P]^\bot H_{PP}\Big)(Du):D^2u\, = \, 0
\eeq
is degenerate elliptic, that is, the tensor map 
\beq \label{4.3}
\A_{\al i \be j}(P)\ := \ H_{P_{\al i}}(P) H_{P_{\be j}}(P) + H(\P)[H_P(P)]_{\al \ga}^\bot H_{P_{\ga i}P_{\be j}}(\P)
\eeq
satisfies the Legendre-Hadamard condition and the symmetry condition
\begin{align}
\A(P):(\eta \ot w) \ot (\eta \ot w)\ \geq \ 0,  \label{4.4}\\
\A(P):(Q \ot R - R \ot Q)\ = \ 0,  \label{4.5}
\end{align}
for all $\eta \in \R^N$, $w\in \R^n$ and $P,Q,R \in \R^{N \by n}$.

\ms

Then, (i) implies (ii). If moreover $H$ is analytic at $0$ and satisfies
\beq \label{4.6} 
\{H=0\}=\{H_P=0\}= 0,\ \ H_{PP}(0) >0\ \& \ H_{PP}:(v\ot w -w \ot v)=0,
\eeq
for $v,w\in \R^n$, $P\in \R^{N\by n}$, then, (ii) implies (i) when either

\noi a) $n\leq 3$,

or

\noi b) $n\geq 4$ and the $q$-th order derivative tensor $H_{P...P}(0) \in \ot^{(q)} (\R^{N\by n})$ is contained in the linear subspace $\msL^q$ which consists of fully symmetric tensors $T$ for which the only non-trivial components are of the form $T_{{\al_1 i}{\al_2 j}{\al_3 k}...{\al_q k}}$, where $\al_m \in \{1,...,N\}$, $i,j,k \in \{1,...,n\}$.

If $n\geq 4$ but $H_{P...P}(0) \not\in \msL^q$, then $H$ has the form \eqref{4.1} up to a fourth order correction: $H(P)=h\big(\frac{1}{2}P^\top P\big) +O(|P|^4)$. If $N\leq n$, it does not follow that $h_p>0$.
\et

\noi In the case that $H(P)$ equals $h\big(\frac{1}{2}P^\top P\big)$, the elliptic system takes the form
\beq \label{4.8}
\A_\infty u \, =\, \Big(Duh_p \ot Duh_p \, +\,  h[Du]^\bot \! \ot h_p\Big) :D^2u \, = \, 0 
\eeq
with $h=h\big(\frac{1}{2}Du^\top Du\big)$.

The extra assumption $H_{P...P}(0) \in \msL^q$ is necessary only in higher dimensions $n\geq 4$. It requires that  $H_{P...P}(0)$ vanishes when more than $3$ of its Latin indices are different to each other. The linear space $\msL^q$ can be descrided as
\begin{align} \label{4.7}
\msL^q\ := \ \Big\{T \in \ot^{(q)} (\R^N\ot \R^n)\ \big| &\ T=T_{{\al_1 i_1}...{\al_q i_q}}e_{\al_1 i_1}\ot ...\ot e_{\al_q i_q} : \nonumber\\
T_{...\al i ...\be j ...} &\ =\ T_{...\be i ...\al j ...} \ =\  T_{...\be j ... \al i ...}\ , \\
&\{i_1,...i_q\} \neq \{i,j,k,...,k\}\ \Longrightarrow \ T=0 
\Big\}.  \nonumber
\end{align}
If $H_{P...P}(0) \not\in \msL^q$, then Hamiltonians with a little more complicated fourth and higher order derivatives also give rise to elliptic systems.

\BPT \ref{th5}. We first prove the implication $(i) \Rightarrow (ii)$. For, assume that the Hamiltonian $H$ has the form \eqref{4.1}. We begin by observing that the symmetry assumption $h_{p_{ij}}=h_{p_{ji}}$ implies that second derivatives of $h$ are fully symmetric in all indices: obviously since $h$ is in $C^2(\mS(\R^n))$ we have $h_{p_{ij}p_{kl}} =h_{p_{kl}p_{ij}} $ 
and also
\beq  \label{4.9}
 h_{p_{ij} p_{kl}} = (h_{p_{kl}})_{p_{ij}} = (h_{p_{lk}})_{ p_{ij} } = h_{ p_{ij} p_{lk} }.
\eeq
By using \eqref{4.9}, we suppress the arguments in the notation of $h$ and calculate
\begin{align} \label{4.10}
H_{P_{\al i}}(P)\ = \ \frac{1}{2}h_{p_{kl}}\Big(\de_{\al \be}\de_{ki}P_{\be l}\ + \ \de_{\al \be}\de_{il}P_{\be k} \Big)  \ =\ P_{\al k}h_{p_{ki}} , 
\end{align}
and
\begin{align} \label{4.11}
H_{P_{\al i} P_{\be j}}(P)\ & = \ \frac{1}{2} h_{p_{ik} p_{lm}}  \Big(\de_{\be \ga}\de_{j l}P_{\ga m}\ + \ \de_{\be \ga}\de_{j m}P_{\ga l}\Big)P_{\al k}  \ +\ h_{p_{ik}}  \de_{\al \be}\de_{k j}   \nonumber\\
&=\  \de_{\al \be} h_{p_{ij}} \ + \ h_{p_{ik} p_{jm}}P_{\al k} P_{\be m}.   
\end{align}
Also, by assumption $[H_P(P)]^\bot=[P]^\bot$. By \eqref{4.10} and \eqref{4.11}, we have
\begin{align}\label{4.13}
\Big(H_P \ot H_P  + H [H_P]^\bot H_{PP}\Big)(P) \ & =\ Ph_p \ot Ph_p \, +\, h[P]^\bot\Big(I \ot h_p\ +\ Ph_{pp}P^\top\Big) \nonumber\\
                               & =\ Ph_p \ot Ph_p \, +\,  h[P]^\bot \ot h_p ,               
\end{align}
where $h=h\big(\frac{1}{2}P^\top P\big)$. Hence, in view of \eqref{4.3}, equation \eqref{4.8} follows. Also, since $h\geq 0$ and $ h_p$, $[P]^\bot$ are positive symmetric, conditions \eqref{4.4} and \eqref{4.5} follow as well:
\begin{align}\label{4.14}
\A(P):(\eta \ot w) \ot (\eta \ot w)  \ &= \ P_{\al k}h_{p_{ki}}\eta_\al w_i P_{\be l}h_{p_{lj}}\eta_\be w_j
 \ + \ h [P]_{\al \be}^\bot \eta_\al \eta_\be h_{p_{ij}}w_i w_j \nonumber\\
& = \ (Ph_p :\eta \ot w)^2 \, +\,  h([P]^\bot  :\eta \ot \eta)( h_p :w\ot w)  \\
& \geq\ 0,                    \nonumber
\end{align}
\begin{align}\label{4.15}
\A(P):(Q \ot R- R\ot Q)  \,  &= \, \pm P_{\al k}h_{p_{ki}}Q_{\al i}  P_{\be l}h_{p_{lj}}R_{\be j}  \ \pm\ h[P]_{\al \be}^\bot Q_{\al i} h_{p_{ij}} R_{\be j} \nonumber\\
 \ &=\, \pm(Ph_p :Q)(Ph_p :R)\  \pm \ h[P]^\bot  : (Q h_p R^\top)  \\
&=\ 0,             \nonumber
\end{align}
for all $\eta \in \R^N$, $w\in \R^n$, $P,Q,R \in \R^{N \by n}$. Hence, $(ii)$ follows. 

Now we assume $(ii)$ and prove the reverse implication. For, suppose $H$ is analytic at $0$ and suppose that  \eqref{4.4} - \eqref{4.7} hold. By \eqref{4.5}, we have
\beq \label{4.16}
\Big(H_P \ot H_P  + H [H_P]^\bot H_{PP}\Big)(P) :(Q \ot R- R\ot Q)\ = \ 0,
\eeq
 for all $P,Q,R \in \R^{N \by n}$. By symmety of $H_P \ot H_P$ and since by \eqref{4.6} we have $H>0$ and $H_P \neq 0$ on $(\R^{N \by n}) \set  \{0\}$, \eqref{4.16} gives
\beq \label{4.17}
 [H_P]^\bot H_{PP} :(Q \ot R- R\ot Q)\ = \ 0.
\eeq
By the identity $[H_P]^\bot = I - [H_P]^\top$ and since $I$, $H_{PP}$ and $[H_P]^\bot$ are symmetric,  for $Q=e_\al \ot e_i$ and $R=e_\be \ot e_j$, \eqref{4.17} gives the commutativity relation
\beq \label{4.18}
 [H_P]_{\al \ga}^\bot H_{P_{\ga i}P_{\be j}}\ = \ H_{P_{\al i}P_{\ga j}}  [H_P]_{\ga \be}^\bot
\eeq
on $(\R^{N \by n}) \set  \{0\}$, that is
\beq \label{4.19}
 [H_P]^\bot H_{PP} \ =\  H_{PP} [H_P]^\bot .
\eeq
We set $A_{\al i \be j} := H_{P_{\al i}P_{\ga j}}(0)$. By assumtpion \eqref{4.6}, we have $A>0$ in $\mS(\R^{N \by n})$. By analyticity of $H$ and since $H(0)=0$ and $H_P(0)=0$, we have
\begin{align}
H(P)\ &=  \ \frac{1}{2}A : P\ot P \ + \ O(|P|^3),  \label{4.20}\\
H_P(P)\ &=  \ A : P \ + \ O(|P|^2),  \label{4.21}\\
H_{PP}(P)\ &=  \ A\ + \ O(|P|), \label{4.22}
\end{align}
as $|P|\ri 0$. Since $A=H_{PP}(0)>0$ and $H_P(0)=0$, the map $H_P$ is a diffeomorphism between open neighbourhods of zero in $\R^{N \by n}$. Hence, there is an $r>0$ such that 
\beq
H_P\ : \ \mB^{Nn}_r  := \{Q \in \R^{N \by n}: |Q|<r\}\, \larrow \, H_P( \mB^{Nn}_r) \, \sub \R^{N \by n}
\eeq
is a diffeomorphism. Hence, there is a $\rho>0$ such that for $0<t<\rho$, $\xi \in \mS^{N-1}$ and $w\in \mS^{n-1}$, there exists a unique $P(t) \in \mB^{Nn}_r$ such that
\beq \label{4.24}
t \, \xi \ot w \ = \ H_P\big(P(t)\big).
\eeq
Moreover, $|P(t)|\ri 0$ as $t\ri 0$. The path $P(\, .\, )$ is the inverse image through $H_P$ of the rank-one line spanned by $\xi \ot w$. By \eqref{4.24}, we have
\beq \label{4.25}
[H_P\big(P(t)\big)]^\top \ = \ [t \, \xi \ot w]^\top\ = \  \xi \ot \xi.
\eeq
By evaluating  \eqref{4.19} at $P(t)$ and using \eqref{4.25} and \eqref{4.22}, we obtain
\beq
(\xi \ot \xi) (A\ + \ o(1))\ = \ (A\ + \ o(1))(\xi \ot \xi) ,
\eeq
at $t\ri 0$. In the limit we get $(\xi \ot \xi) A\ = \ A(\xi \ot \xi )$, that is
\beq
\xi_\al \xi_\ka A_{\ka i \be j}\ = \  A_{\al i \ka j} \xi_\ka \xi_\be,
\eeq
for all $i,j \in \{1,...,n\}$, $\al ,\be \in \{1,...,N\}$. By the symmetry condition in assumption \eqref{4.6}, for all $i,j$ fixed the matrix $ A_{\al i \be j}$ commutes with all $1$-dimensional projections of $\R^N$. Hence, it is simultaneously diagonalisable with them and as such a multiple of the identity. Thus, there is a symmetric matrix $\hat{A}_{ij}$ such that 
\beq \label{4.28}
 A_{\al i \be j}\ = \ \hat{A}_{ij} \de_{\al \be}.
\eeq
Consequently, $A:P\ot P = \hat{A} : P^\top P$. We now set $B_{\al i \be j \ga k}:= H_{P_{\al i}P_{\be j}P_{\ga k}}(0)$. Then, by \eqref{4.28}, equations \eqref{4.21} and \eqref{4.22} become
\begin{align}
H_{P}(P)\ &=\ P\hat{A} \ + \ O(|P|^2), \label{4.29}\\
H_{PP}(P)\ &=\ I \ot \hat{A} \ + \frac{1}{2}B:P \ + \ O(|P|^2),\label{4.30}
\end{align}
and hence by \eqref{4.29} and \eqref{4.24} we get 
\beq \label{4.31}
t\, \xi \ot w\ =  \ P(t)\hat{A}\ +\ O(|P(t)|^2).
\eeq
Since $A>0$ in $\mS(\R^{N \by n})$, we have $\hat{A}>0$  in $\mS(\R^n)$ as well. Thus, for $0<t<\rho$, we have
\beq \label{4.32}
\frac{P(t)}{|P(t)|}\ +\ O(|P(t)|)\ = \ \frac{t}{|P(t)|}\xi \ot \big((\hat{A}^{-1})^\top w\big).
\eeq
As $t\ri 0$, we have $|P(t)|\ri 0$ and by compactness along an infinitessimal sequence $t_m \ri 0$ there exists a $\bar{P}$ with $|\bar{P}|=1$ such that $P(t_m)/|P(t_m)| \ri \bar{P}$. By passing to the limit in \eqref{4.32} as $m\ri \infty$ along $\{t_m\}$, we obtain that the limit of ${t_m}/{|P(t_m)|}$ exists
and
\beq
\label{4.33}
\lim_{m \ri \infty}\frac{P(t_m)}{|P(t_m)|}\ =\ \bar{P}\ = \ \xi \ot \left[\left(\lim_{m \ri \infty}\frac{t_m}{|P(t_m)|}\right)(\hat{A}^{-1})^\top w \right].
\eeq
Since $\hat{A}^{-1}>0$ and $|\bar{P}|=1$, for any $v \in \mS^{n-1}$, there is a $w \in \mS^{n-1}$ such that \eqref{4.33} becomes
\beq
\label{4.34}
\lim_{m \ri \infty}\frac{P(t_m)}{|P(t_m)|}\ =\ \bar{P}\ = \ \xi \ot v.
\eeq
By \eqref{4.19}, \eqref{4.25}, \eqref{4.29} \eqref{4.30}, we have
\begin{align} \label{4.35}
\xi \ot \xi\Big( I \ot \hat{A} \ + \ &\frac{1}{2}B:P(t) \ + \ O(|P(t)|^2) \Big) \nonumber\\
& = \ \Big( I \ot \hat{A} \ +\ \frac{1}{2}B:P(t) \ + \ O(|P(t)|^2) \Big)\xi \ot \xi.
\end{align}
By cancelling the commutative term $\xi \ot \xi( I \ot \hat{A} )$, \eqref{4.35} gives
\begin{align} \label{4.36}
\xi \ot \xi\left( B:\frac{P(t)}{|P(t)|} \ + \ O(|P(t)|) \right) 
 = \ \left( B:\frac{P(t)}{|P(t)|} \ + \ O(|P(t)|) \right) \xi \ot \xi.
\end{align}
By passing to the limit in \eqref{4.36} along $t_m \ri 0$, in view of \eqref{4.34} we obtain
\beq \label{4.37}
\xi \ot \xi (B : \xi \ot v)\ = \  (B : \xi \ot v)\xi \ot \xi,
\eeq
for all $\xi \in \mS^{N-1}$, $v \in \mS^{n-1}$. Hence, \eqref{4.37} for $v=e_k$ says
\beq \label{4.40}
\xi_\al (B_{\be i \la j \mu k} \, \xi_\mu  \xi_\la ) \ =\ (B_{\al i \la j \mu k}\, \xi_\mu \xi_\la)   \xi_\be.
\eeq
By \eqref{4.40}, $B: \xi \ot \xi$ is proportional to $\xi$; hence, there is a map $\hat{B} : \R^N \larrow \ot^{(3)}\R^n$ such that $B: \xi \ot \xi = \hat{B}(\xi) \ot \xi$, or
\beq \label{4.41}
B_{\al i \la j \mu k} \, \xi_\mu  \xi_\la \ = \ \hat{B}_{ijk}(\xi)\, \xi_\al. 
\eeq
By assumption \eqref{4.6} and induction, all second and higher order derivatives are fully symmetric in all their indices. Hence, we may fix $i,j,k \in \{1,...,n\}$ and suppress the dependence in them to obtain $B_{\al \ka \la } \, \xi_\ka  \xi_\la \ = \ \hat{B}(\xi)\, \xi_\al$ with $\hat{B} \in C^\infty(\R^N \set \{0\})$. The idea now is to differentiate in order to cancel both $\xi$'s contracted with $B$ and then contract again with a vector which annihilates $\xi$ from the right hand side. For, by differentiating we get
\beq  \label{4.42}
D_\be\hat{B}(\xi)\xi_\al \ = \ - \hat{B}(\xi)\de_{\al \be}\ +  \ 2 B_{\al \be \ga} \xi_\ga.
\eeq
By \eqref{4.42}, we obtain that $D\hat{B}(\xi)\ot \xi$ is symmetric. Hence, we get that $D\hat{B}(\xi)\ot \xi=\xi \ot D\hat{B}(\xi)$ and hence there exists $\bar{B} \in C^\infty (\R^N \set \{0\})$, such that $D\hat{B}(\xi) =\bar{B}(\xi)\xi$. Thus, \eqref{4.42} gives
\beq \label{4.44}
\bar{B}(\xi)\xi \ot \xi \ +\hat{B}(\xi)\, I\ = \ 2 B :\xi.
\eeq
By differentiating the expression $D\hat{B}(\xi) =\bar{B}(\xi)\xi$, we get
\beq \label{4.45}
D\bar{B}(\xi)\ot \xi \ =\ D^2\hat{B}(\xi)\ - \ \bar{B} (\xi)\, I.
\eeq
By \eqref{4.45}, we obtain that $D\bar{B}(\xi)\ot \xi$ is symmetric too. Hence, there exists $\check{B} \in C^\infty (\R^N \set \{0\})$, such that $D\bar{B}(\xi) =\check{B}(\xi) \xi$ and hence \eqref{4.45} gives
\beq \label{4.46}
D^2\hat{B}(\xi) \ = \ \check{B}(\xi)\xi \ot \xi\ + \ \bar{B} (\xi)\, I.
\eeq
By differentiating \eqref{4.42} again and inserting \eqref{4.46} we get
\begin{align} \label{4.47}
2B_{\al \be \ga}\ &= \ D^2_{\be \ga}\hat{B}(\xi)\xi_\al\ + \ D_\be\hat{B}(\xi)\de_{\al \ga}\ + \ D_\ga \hat{B}(\xi) \de_{\al \be}\\
& = \ \check{B}(\xi)\xi_\al \xi_\be \xi_\ga \ + \ \bar{B}(\xi) \big(\xi_\al \de_{\be \ga}\ + \ \xi_\be \de_{\al \ga}\ +\ \xi_\ga \de_{\be \al}\big), \nonumber
\end{align}
for all $\xi \in \R^N \set \{0\}$. Since $N\geq 2$, for each $\eta \in\R^N$ we can choose a nonzero $\xi$ normal to $\eta$. Hence, by triple contraction in \eqref{4.47} we obtain
\begin{align} \label{4.48}
B: \eta \ot \eta \ot \eta \  = \ \frac{1}{2} \Big[ \check{B} (\xi) ( \xi^\top \eta )^2 \ + \ \bar{B} (\xi) | \eta |^2 \Big]( \xi^\top \eta) \  =\ 0.
\end{align}
Hence, by full symmetry in all indices we obtain $H_{P_{\al i} P_{\be j} P_{\ga k}}(0)= B_{{\al i} {\be j} {\ga k}}=0$ and consequently third order derivatives vanish. We now set
\beq \label{4.50}
C_{{\al i} {\be j} {\ga k} {\de l}} \ := \ H_{P_{\al i} P_{\be j} P_{\ga k} P_{\de l}}(0)
\eeq
and then for $0<t<\rho$, \eqref{4.30} and \eqref{4.36} become
\begin{align}
H_{PP}(P(t))\ =& \ I  \ot \hat{A} \ +\ \frac{1}{3!}C:P(t)\ot P(t)\ + \ O(|P(t)|^3),\\
\xi \ot \xi & \left( C:\frac{P(t)}{|P(t)|}\ot\frac{P(t)}{|P(t)|} \ + \ O(|P(t)|) \right) 
\\
&\ \ \ \ = \ \left( C:\frac{P(t)}{|P(t)|} \ot \frac{P(t)}{|P(t)|} \ + \ O(|P(t)|) \right) \xi \ot \xi. \nonumber
\end{align}
By setting $t=t_m$ and letting $m\ri \infty$, in view of \eqref{4.34}, we get
\begin{align} \label{4.52}
\xi \ot \xi \Big[C:(\xi \ot v) \ot (\xi \ot v)\Big] \ = \ \Big[C:(\xi \ot v) \ot (\xi \ot v)\Big]\xi \ot \xi ,
\end{align}   
for all $\xi \in \R^N$, $v \in \R^n$. Hence, for $v=e_k$,
\begin{align} \label{4.53}
\xi_\al  \Big[C_{\be i \ka j \la k \mu k}\, \xi_\ka \xi_\la \xi_\mu\Big] \ = \ \Big[C_{\al i \ka j \la k \mu k}\, \xi_\ka \xi_\la \xi_\mu\Big]  \xi_\be .
\end{align}   
By \eqref{4.53}, there exists a map $\hat{C} \, :\, \R^N \set \{0\} \larrow \ot^{(4)}\R^n$ with $\hat{C}_{ijkk} \in C^\infty(\R^N \set \{0\})$ such that
\beq
C_{\al i \ka j \la k \mu k}\, \xi_\ka \xi_\la \xi_\mu \ = \ \hat{C}_{ijkk}(\xi)\xi_\al.
\eeq
By fixing again the indices $i,j,k$, dropping them and arguing exactly as we did before for $B_{\al \be \ga}$, there exist $\bar{C},\check{C} \in C^\infty(\R^N \set \{0\})$ such that
\begin{align} \label{4.55}
3! C_{\al \be \ga \de}\xi_\de\  = \ \check{C}(\xi)\xi_\al \xi_\be \xi_\ga \ + \ \bar{C}(\xi) \big(\xi_\al \de_{\be \ga}\ + \ \xi_\be \de_{\al \ga}\ +\ \xi_\ga \de_{\be \al}\big).
\end{align}
By differentiating \eqref{4.55}, we get
\begin{align} \label{4.56}
3!\, C_{\al \be \ga \de}\ -& \ \bar{C}(\xi) \big(\de_{\al \be} \de_{\ga \de}\ + \ \de_{\ga \be} \de_{\al \de}\ +\ \de_{\de \be} \de_{\ga \al} \big) \nonumber\\
 =& \ \check{C}(\xi)  \big( \xi_\al \xi_\be \de_{\ga \de} \ + \  \xi_\be \xi_\ga \de_{\al \de}\ +\ \xi_\ga \xi_\al \de_{\be \de} \big)\\
& \ + D_\de \check{C}(\xi)\Big[ \xi_\al \xi_\be \xi_\ga \ + \ \big(\xi_\al \de_{\be \ga}\ + \ \xi_\be \de_{\al \ga}\ +\ \xi_\ga \de_{\be \al}\big)\Big]. \nonumber
\end{align}
Fix $\eta \in \R^N$. Since $N\geq 2$, there exists $\xi \perp \eta$, $\xi \neq 0$. Then, \eqref{4.56} gives
\begin{align} \label{4.57}
\left[C_{\al \be \ga \de}\ - \ \frac{\bar{C}(\xi)}{3!} \big(\de_{\al \be} \de_{\ga \de}\ + \ \de_{\ga \be} \de_{\al \de}\ +\ \de_{\de \be} \de_{\ga \al} \big) \right]
 \eta_\al \eta_\be \eta_\ga \eta_\de \ &= \ O(|\eta^\top \xi|) \nonumber\\
&=\ 0.
\end{align}
By \eqref{4.57}, the function $\bar{C}$ is constant and moreover for all $i,j,k$,
\beq \label{4.58}
C_{\al i \be j  \ga k \de k}\ = \ \frac{\bar{C}_{ijkk}}{3!} \big(\de_{\al \be} \de_{\ga \de}\ + \ \de_{\ga \be} \de_{\al \de}\ +\ \de_{\de \be} \de_{\ga \al} \big).
\eeq
If either $n\leq 3$ or $n\geq 4$ but $H_{PPPP}(0) \in \msL^4$, where $\msL^4$ is given by \eqref{4.7}, then in view of \eqref{4.50}, the tensor $C_{\al i \be j  \ga k \de l}$ has no more than $3$ different indices $i,j,k,l$ for which it is non-zero. Hence, by full symmetry in all indices, \eqref{4.58} completely determines $H_{PPPP}(0)$ and gives
\begin{align}
H_{PPPP}(0):\ot^{(4)}P \ = \ \frac{1}{2}\bar{C}_{ijkl}P_{\al i}P_{\al j}P_{\be k}P_{\be l} \ = \  \frac{\bar{C}}{2}: (P^\top P)\ot (P^\top P).
\end{align}
Now we iterate the above arguments. The analog of \eqref{4.52} after blowing up along $t_m$ for $q$-th order derivatives is
\beq \label{4.60}
\xi \ot \xi \Big[H_{P...P}(0):\ot^{(q-2)} (\xi \ot v)\Big] \ = \ \Big[ H_{P...P}(0) :\ot^{(q-2)} (\xi \ot v) \Big]\xi \ot \xi,
\eeq
for all $\xi \in \R^N$, $v\in \R^n$. When $H_{P...P}(0) \in \msL^q$, the only components of the tensor $H_{P_{\al_1 i_1}...P_{\al_q i_q}}(0)$ which may not vanish are of the form 
\beq
H_{P_{\al_1 i}P_{\al_2 j}P_{\al_3 k}...P_{\al_q k}}(0),
\eeq
where $i,j,k \in \{1,...,n\}$ and $\al_1,...,\al_q \in \{1,...,N\}$. Hence, \eqref{4.60}, completely determines $H_{P...P}(0)$. By induction, all odd order derivatives of $H$ vanish and all even order derivatives depend on $P$ via $P^\top P$: we have
\beq
\underset{\text{q-th order}}{\underbrace{H_{P...P}(0)}}:\ot^{(q)}P\ = \ 
\left\{
\begin{array}{l}
C_q :\ot^{(q/2)} P^\top P, \ \ q \in 2\N, \ms\\
0, \hspace{70pt} q \in 2\N +1,
\end{array}
\right.
\eeq
for certain tensors $C_q \in \ot^{(q/2)}\R^n$. Hence, by defining $h : \R^{n \by n} \larrow \R$ by 
\beq
h(p)\ :=\ \sum_{m=1}^\infty 2^m C_{2m} :\ot^{(m)}p,
\eeq
we obtain
\beq \label{4.64}
H(P) \ =\ h \Big(\frac{1}{2}P^\top P\Big).
\eeq
Hence, $h\geq 0$ with $h\in C^\infty\big(\mS(\R^n)^+\big)$ and also $h_p = h_p^\top$. Moreover, by assumption and \eqref{4.64}  we have $[P]^\bot=[H_P(P)]^\bot=[Ph_p(\frac{1}{2}P^\top\! P)]^\bot$.

If finally $H_{P...P}(0) \not\in \msL^q$, then $H$ has the form \eqref{4.64} up to a correction of order $O(|P|^4)$. This follows by decomposing each $H_{P...P}(0)$ to the sum of a term in $\msL^q$ and a term in the orthogonal complement of  $\msL^q$. The $O(|P|^4)$ function arises from the series consisting of the forth and higher order parts of $H_{P...P}(0):\ot^{(q)}P$ in the orthogonal complements. The theorem follows.                                         \qed

\section{The $1$-dimensional case of ODE system with dependence on all arguments.} \label{section5}

\subsection{Formal derivation of the general ODE System.}  \label{subsection5.1} Let $H$  be a non-negative Hamiltonian in $C^2(\R \by \R^N \by \R^N)$, where $N\geq 2$ and we denote the arguments of $H$ by $H(x,\eta,P)$. Consider the integral functional 
\beq \label{5.1}
E_m(u,I)\ := \ \int_I \big(H(x,u(x),u'(x))\big)^mdx,
\eeq
where $m\geq 2$ and $u : I \sub \R \larrow \R^N$. The Euler-Lagrange equation of functional \eqref{5.1} is the ODE system
\beq \label{5.2}
\Big(H^{m-1}(\cdot , u,u')H_P(\cdot , u,u')\Big)' \ = \ H^{m-1}(\cdot , u,u')H_\eta (\cdot , u,u')
\eeq
which after expansion and normalisation gives
\beq \label{5.3}
\big(H(\cdot , u,u')\big)' H_P(\cdot , u,u') \ + \ \frac{H(\cdot , u,u')}{m-1}\Big( \big(H_P(\cdot , u,u')\big)' - H_\eta (\cdot , u,u')\Big) \ = \ 0,
\eeq
on $I\sub \R$. We define the following projections of $\R^N$:
\begin{align}
[H(x,\eta,P)]^\top \ & := \ \sgn\big(H_P(x,\eta,P)\big) \ot  \sgn\big(H_P(x,\eta,P)\big), \label{5.4}\\
[H(x,\eta,P)]^\bot \ & := \ I \ - \ [H(x,\eta,P)]^\top. \label{5.5}
\end{align}
Then, by employing \eqref{5.4} and \eqref{5.5} to expand the term in bracket of \eqref{5.3}, we obtain
\begin{align}
\big(H(\cdot , u,u')\big)' H_P(\cdot , u,u') \ +& \ \frac{H(\cdot , u,u')}{m-1}[H_P(\cdot , u,u')]^\top\Big(\big(H_P(\cdot , u,u')\big)' 
- H_\eta (\cdot , u,u')\Big)\label{5.6}\\
= -& \ \frac{H(\cdot , u,u')}{m-1}[H_P(\cdot , u,u')]^\bot\Big(\big(H_P(\cdot , u,u')\big)' 
- H_\eta (\cdot , u,u')\Big) .
\nonumber
\end{align}
By perpendicularity of the orthogonal projections \eqref{5.4} and \eqref{5.5}, the left and right hand sides of \eqref{5.6} are normal to each other. Hence, they both vanish. By re-normalising the right hand side and rearranging, we get
\begin{align}
\big(H(\cdot , u,u')\big)' H_P(\cdot , u,u') \ +& \ H(\cdot , u,u')[H_P(\cdot , u,u')]^\bot\Big(\big(H_P(\cdot , u,u')\big)' - H_\eta (\cdot , u,u')\Big) \label{5.7} \nonumber\\
= -& \ \frac{H(\cdot , u,u')}{m-1}[H_P(\cdot , u,u')]^\top\Big(\big(H_P(\cdot , u,u')\big)' 
- H_\eta (\cdot , u,u')\Big) .
\end{align}
As $m\ri \infty$, we obtain the complete system of fundamental ODEs for a general Hamiltonian with dependence on all the arguments
\begin{align}
\big(H(\cdot , u,u'&)\big)' H_P(\cdot , u,u') \ +\  H(\cdot, u,u') \cdot \nonumber\\ 
&\cdot[H_P(\cdot, u,u')]^\bot\, \Big(\big(H_P(\cdot, u,u')\big)' - H_\eta (\cdot, u,u')\Big)  =\ 0 \label{5.8},
\end{align}
whose solutions are curves $u: I \sub \R \larrow \R^N$.

\ms

\subsection{Degenerate elliptic ODE systems.} \label{subsection5.2} We begin by observing that the Ellipticity Classification Theorem \ref{th5} readily extends to the case of $H(x,\eta,P)$ with dependence on all arguments; the form \eqref{4.1} of the Hamiltonian modifies to 
\beq \label{5.9} 
H(x,\eta,P)\ = \ h\Big(x,\eta,\frac{1}{2}P^\top P\Big)
\eeq
and the PDE systems \eqref{4.2} and \eqref{4.8} modify by the appearance of first and lower order tems. In the case of ODEs where $n=1$, the ``geometric" Hamiltonians of the form \eqref{5.9} become the \emph{radially symmetric} ones:
\beq \label{5.10} 
H(x,\eta,P)\ = \ h\Big(x,\eta,\frac{1}{2}|P|^2\Big),
\eeq
where $h \in C^2\big(\R \by \R^N \by [0,\infty)\big)$ and the degenerate elliptic ODE system takes a particularly important and simple form. We note that when we have lower order terms, the Hamiltonian
\[
H(x,\eta,P)\ = \ h\Big(x,\eta,\frac{1}{2}|P-V(x,\eta)|^2\Big)
\]
also leads to degenerate elliptic system, and this is important elsewhere \cite{K9}. However, for simplicity herein we choose $V\equiv 0$. In the case of $\De_\infty$, we have $h(x,\eta,p)=p$.  We now derive the ODEs in the elliptic case.

Suppose $h \in C^2\big(\R \by \R^N \by [0,\infty)\big)$ with arguments denoted by $h(x,\eta,p)$ and define $H \in C^2(\R \by \R^N \by \R^N)$ by means of \eqref{5.10}. We henceforth assume
\beq \label{5.11}
\big\{h_p(x, \eta,\cdot )=0\big\}\, \sub \, \{0\} \, =\, \big\{h(x, \eta,\cdot )=0\big\},
\eeq
for all $(x,\eta)\in \R^{1+N}$. Assumption \eqref{5.11} is natural and will make the normal coefficient $H[H_P]^\bot$ of \eqref{5.8} continuous. By using \eqref{5.10} and supressing arguments, we calculate derivatives:
\begin{align} 
H_P\ &= \ h_p P,\ \ \ \ \ \ \ \ H_{PP} \ = \ h_{pp}P \ot P \ +\ h_p I, \ \ H_\eta \ =  \ h_\eta, \label{5.12}\\
H_{P\eta} \ &=\ P\ot h_{p\eta}, \ \  H_{Px}\ = \ h_{px}P, \hspace{65pt} H_x \ = \  h_x.\label{5.13}
\end{align}
By expanding derivatives in \eqref{5.8} and using \eqref{5.10}, \eqref{5.12} and \eqref{5.13}, we get
\begin{align} 
(h_p)^2(u' \ot u' )u''\ &+\ h_p (u' \ot h_\eta )u' \ + \ h_xh_p u' \nonumber  \\
+\, h[h_p u']^\bot\Big(& h_{pp} (u' \ot u')u'' \ +\ (u' \ot h_{p\eta})u' \ \label{5.15}\\
& +\ h_{px}u'\  +\ h_p u'' \ -\ h_\eta\Big)\ = \  0,\nonumber
\end{align}
where $h=h\big(\cdot ,u,\frac{1}{2}|u'|^2\big)$. By assumption, \eqref{5.11}, we have $\{h_p u' =0\} =\{u'=0\}= \{h=0\}$. Hence, we obtain that $[h_p u']^\bot = [u']^\bot$. On $\{u'\neq 0\}$, we multiply the normal term of \eqref{5.15} by $\frac{|u'|^2 h_p}{h}$ to obtain
\begin{align}  \label{5.16}
(h_p)^2(u' \ot u' )u''\ &+\ h_p \Big((u' \ot u')h_\eta \ + \ h_x u' \Big)\nonumber  \\
&+\ |u'|^2 h_p[u']^\bot\Big( h_p u'' \ - \ h_\eta \Big) \ = \ 0. 
\end{align}
Hence, by using the identity $|u'|^2 I = u' \ot u' + |u'|^2 [u']^\bot $, \eqref{5.16} gives
\beq \label{5.17}
(h_p)^2 |u'|^2 u''\ -\ h_p \left(|u'|^2 \Big( I - 2 \frac{u'}{|u'|} \ot  \frac{u'}{|u'|} \Big) h_\eta \ - \ h_x u' \right) \ = \ 0.
\eeq
By introducing the \emph{reflection operator $\bR_\xi : \R^N \larrow \R^N$ with respect to the hyperplane $[\xi]^\bot$, $\xi \in \R^N \set \{0\}$}, given by
\beq \label{5.18}
\bR_\xi \ := \ I\ -\ 2 \frac{\xi}{|\xi|} \ot  \frac{\xi}{|\xi|},
\eeq
the ODE system \eqref{5.17} becomes
\beq \label{5.19}
A_\infty u \ := \ |u'|^2\Big(h_p u''\, - \, \bR_{u'}h_\eta\Big) \ + \ h_x u'\ = \  0 ,
\eeq
where $h=h\big(\cdot ,u,\frac{1}{2}|u'|^2\big)$. In view of \eqref{5.11}, the systems \eqref{5.19} and \eqref{5.8} are equivalent on $\{u'=0\}$ as well.  The system \eqref{5.19} comprises the \emph{degenerate elliptic ODE system}.

\begin{remark}
We observe that in the special case where $h=h\big(\frac{1}{2}|u'|^2\big)$ and $h_\eta \equiv 0$, $h_x \equiv 0$, solutions of \eqref{5.19} trivialize to \emph{affine} and actually \eqref{5.19} is equivalent to $\De_\infty$. In the special case where $h=h\big(\cdot,\frac{1}{2}|u'|^2\big)$ and $h_\eta \equiv 0$, solutions of \eqref{5.19} become essentially scalar with \emph{affine rank-one range}, that is $u(\R)$ is contained in an affine line of $\R^N$ since $u''$ becomes proportional to $u'$ and \eqref{5.19} becomes essentially scalar. Consequently, \eqref{5.19} is most interesting when $h\big(x,u(x),\frac{1}{2}|u'(x)|^2\big)$ depends on $u(x)$ and hence $h_\eta \not\equiv 0$. In this case the reflection operator $\bR_{u'}$ with respect to the normal hyperplane $[u']^\bot$ is \emph{discontinuous on $\{u'=0\}$ at critical points of $u$}, but the product $|u'|^2\bR_{u'}$ is continuous. However, in all cases the system is always degenerate.
\end{remark}

\subsection{The initial value problem for the elliptic ODE systems.} \label{subsection5.3} 

In this subsection we solve the initial value problem for ODE system \eqref{5.19} and consider some regularity questions.

\begin{theorem}[The initial value problem for the ODE system] \label{th6} Suppose that $h \in C^2\big(\R \by \R^N \by [0,\infty)\big)$ satisfies $h,h_p \geq 0$ and also \eqref{5.11} and consider the following problem
\beq \label{5.24}
\left\{
\begin{array}{c}
A_\infty u  \ = \ |u'|^2\Big(h_p u''\, - \, \bR_{u'}h_\eta\Big) \ + \ h_x u'\ = \  0 , \ms\\
u(x_0) = u_0\ ,\ \ \ u'(x_0)=v_0\ , \ \ x_0 \in \R.
\end{array}
\right.
\eeq
Then: 

\noi (i) For any \emph{non-critical} initial conditions $(u_0,v_0)\in \R^N \by (\R^N \set \{0\})$, there exists a unique maximal smooth solution $u : (x_0-r,x_0+r)\larrow \R^N$ for some $r>0$ which solves \eqref{5.24} and satisfies $|u'|>0$.

\ms

\noi(ii) For any \emph{critical} initial condition $(u_0,0)\in \R^N \by \{0\}$, there exists at least one solution to \eqref{5.24}, one of them being the constant one $u\equiv u_0$.

\ms

\noi (iii) If
\beq \label{5.25}
h_\eta(x,\eta,0) \neq 0 \ \ \ \text{ and }\ \ \ h_x(x,\eta,p) = O(p) \ \text{ as }\  p\ri 0, 
\eeq
then bounded maximal solutions of\eqref{5.24} starting (in positive time) from non-critical data, either are defined on $[x_0,\infty)$ being smooth and satisfying $|u'|>0$, or they reach a critical point $u'=0$ and form a discontinuity in $u''$ in finite time.

\ms

\noi (iv) If
\beq \label{5.25a}
c\, \leq\, h_p\, \leq\, \frac{1}{c} \text{ \ \  for $c>0$,\ \ and }\ \ \ h_x(x,\eta,p) = O(p) \ \text{ as }\  p\ri 0, 
\eeq
then bounded maximal solutions of \eqref{5.24} either are globally smooth or can be extended past singularities as $W^{2,\infty}_{loc}(\R)^N$ strong solutions which satisfy \eqref{5.19} everywhere and are eventually constant.
\end{theorem}

The interpretation of $W^{2,\infty}_{loc}(\R)^N$ solutions to \eqref{5.24} as strong everywhere  solutions is the same as in Aronsson \cite{A1, A2, A5}: at critical points of $u$ whereon $u''$ may not exist but is essentially bounded in a neighbourhod of $\{u'=0\}$, the coefficient $|u'|^2$ vanishes.

\begin{example} The solution of problem \eqref{5.24} is generally \emph{non-unique} for critical initial conditions. Choose $h(x,\eta,p):=\frac{1}{2}|\eta|^2+p$. Then, \eqref{5.19} takes the form
\beq \label{5.27}
|u'|^2\Big(u'' \, -\, \bR_{u'}u\Big)\ = \ 0
\eeq
and the Hamiltonian is $H(u,u')=\frac{1}{2}(|u|^2+|u'|^2)$. In view of example $3$ in Aronsson's paper \cite{A1}, for essentially scalar solutions $u=\xi v$ where $\xi \in \mS^{N-1}$ and $v:\R\larrow \R$, \eqref{5.27} takes the form $|v'|^2\big(v'' +v\big)\xi=0$. Hence, for initial conditions $u\big(-\frac{\pi}{2}\big)=-e_1$, $u'\big(-\frac{\pi}{2}\big)=0$, \eqref{5.20} admits the solutions $u_1(x)=e_1\sin x$ and $u_2(x) = -e_1$.
\end{example}

\noi The non-uniqueness for critical data owes to that \eqref{5.19} is an $1$-dimensional degenerate elliptic system and the initial value problem is not well-posed for it.

\BPT \ref{th6}. All assertions follow directly by considering the following dynamical formulation of the ODE \eqref{5.19}. For, we write the $N$-dimensional second order degenerate implicit system \eqref{5.19} as a $2N$-dimensional first order explicit system for a vector field defined off an $N$-dimensional ``slice'' of $\R^{2N}$. For $U=(u,v)^\top \in \R^{2N}$, we set
\beq \label{5.20}
 U(x)\ := \ (u(x),u'(x))^\top , \ \ \ \ U  : I \sub \R \larrow \R^{2N},
\eeq
\beq \label{5.21}
F(x,U)\ := \ \left[
                   \begin{array}{c}v\\
\frac{1}{h_p\big(x,u,\frac{1}{2}|v|^2\big)}\left( \bR_v h_\eta \big(x,u,\frac{1}{2}|v|^2\big)  \ - \ \frac{h_x\big(x,u,\frac{1}{2}|v|^2\big)}{|v|^2}v\right)
                       \end{array}
                       \right],
\eeq
where
\beq
F\ : \ \R \by \R^N \by \big(\R^N \set \{0\}\big)\ \larrow \ \R^{2N}.
\eeq
Then, in view of \eqref{5.20} and \eqref{5.21}, ODE system \eqref{5.19} can be written as 
\beq \label{5.23}
\dot{U}(x)\ = \ F\big(x,U(x)\big)\ , \ \ \ \ U  : I \sub \R \larrow \R^{2N}.
\eeq
We now merely observe that the equation
\beq
u'' \ =  \ \frac{1}{h_p\big(\cdot,u,\frac{1}{2}|u'|^2\big)}\left( \bR_{u'} h_\eta \Big(\cdot,u,\frac{1}{2}|u'|^2\Big)  \ - \ \frac{h_x\big(\cdot,u,\frac{1}{2}|u'|^2\big)}{|u'|^2}u'\right)
\eeq
which follows by \eqref{5.23}, implies that under assumption \eqref{5.25} the first term in the bracket becomes discontinuous at critical points of $u$, while the second one vanishes. Solutions extend past critical points where $u''$ ``jumps'' by constant solutions. \qed

\ms

\section{Rigidity of radial $2$-dimensional solutions.} \label{section6}

In this section we study a class of special solutions of the $\infty$-Laplacian, that of smooth $\infty$-Harmonic maps $u : \R^2 \larrow \R^N$, $N\geq 2$ of the form $u=\rho^k f(k\theta)$ in polar coordinates $(\rho,\theta)$. Here $k>0$ is a parameter and $f : \R \larrow \R^N$ is a curve in $\R^N$. It follows that such solutions are very rigid, because  if $k\neq 1$ they are essentially scalar and if $k=1$ they always have affine image. The result here is 

\begin{proposition}[Rigidity of radial 2D $\infty$-Harmonic maps] \label{Pr1} Let $u : \R^2 \larrow \R^N$ be an $\infty$-Harmonic map of the form $u=\rho^k f(k\theta)$ in polar coordinates $(\rho,\theta) \in \R^2$, $k>0$, $f\in C^\infty(\R)^N$, $N\geq 2$. Then, $f$ solves the ODE systems
\begin{align}
f' \ot f' \big( f'' \, + \, f\big) \ +\ \frac{k-1}{k} \big( |f'|^2\, +\, |f|^2 \big) f\ = \ 0, \label{3.21}\\ 
[(f' , f )]^\bot f'' \ = \ 0 .\label{3.22}
\end{align}
Moreover:

\noi (i) If $k\neq 1$, then all solutions have constant rank one, the image $u(\R^2)$ is contained into a line passing through the origin and $f$ can be represented as $f(\theta)=\xi g(\theta)$ for some  $\xi \in \mS^{N-1}$ and $g\in C^\infty(\R)$. 

\noi (ii) If $k=1$, then all solutions have rank at most two and the image $u(\R^2)$ is contained into a $2$-plane of $\R^N$ passing through the origin. On this plane $f$ can be represented as
\beq
f(\theta)\ = \ c \cos B(\theta)\, \big(\cos A(\theta), \sin A (\theta)\big)^\top,
\eeq
where  $c\in \R$ and $A,B \in C^\infty(\R)$ satisfy $|B'|^2 + |A'|^2 \cot^2 B = 1$ and $0<B\leq\frac{\pi}{2}$.
\end{proposition}

\BPP \ref{Pr1}. The derivation of the ``tangential part" \eqref{3.21} of $\De_\infty$ is entirely analogous to Aronsson's derivation of its scalar counterpart in the paper \cite{A6}, p.\ 138. Hence, it suffices to outline the derivation of the ``normal part" \eqref{3.22}. Since for all $\al \in \{1,...,N\}$ we have $u_\al =\rho^k f_\al(k\theta)$, we obtain
\begin{align} \label{3.25}
\left[
\begin{array}{c}
D_x u_\al\\
D_y u_\al
\end{array}
\right]\ &= 
\ \left[
\begin{array}{cc}
\cos \theta & -\sin \theta\\
\sin \theta & \cos \theta
\end{array}
\right]
\left[
\begin{array}{c}
D_\rho u_\al\\
\frac{1}{\rho} D_\theta u_\al
\end{array}
\right]     \nonumber\\
&= 
\ \left[
\begin{array}{cc}
\cos \theta & -\sin \theta\\
\sin \theta & \cos \theta
\end{array}
\right]
\left[
\begin{array}{c}
k\rho^{k-1} f_\al\\
k\rho^{k-1} f'_\al
\end{array}
\right].
\end{align}
 Let $O(\theta)$ denote the rotation-by-$\theta$ appearing in \eqref{3.25}. By recalling that $(f,f')$ is a matrix-valued curve $\R \larrow \R^N \ot \R^2$, we may write
\beq
Du \ = \ k\rho^{k-1}(f,f')\, O(\theta)^\top.
\eeq
Hence, since $ O(\theta)^\top =  O(\theta)^{-1}$ we have
\begin{align}
N(Du^\top)\ = \ \big\{\eta \in \R^N\ : \ \eta^\top (f,f')\, O(\theta)^\top=0 \big\} \ =\ N\big((f,f')^\top\big). 
\end{align}
and consequently $[Du]^\bot=[(f,f')]^\bot$. Moreover,
\begin{align} \label{3.26}
[Du]^\bot\De u \ &= [(f,f')]^\bot \left(\frac{1}{\rho} D_\rho u \, +\,  D^2_{\rho \rho} u\, +\, \frac{1}{\rho^2} D^2_{\theta \theta} u \right) \nonumber\\
         &= [(f,f')]^\bot \Big(k\rho^{k-2}f\, +\, k(k-1)\rho^{k-2}f\, +\,  k^2\rho^{k-2}f''\Big)   \\
          &= k^2\rho^{k-2} [(f,f')]^\bot f''. \nonumber
\end{align}
By Corollary \ref{Cor3}, we may require $|Du|>0$ and hence \eqref{3.22} follows by \eqref{3.8} and \eqref{3.26}. Now, for $(i)$ we have that if $k\neq 1$ then on $\{|f|>0\}$ \eqref{3.21} gives
\beq
-\frac{k\big( f'' \, + \, f\big)^\top f'}{(k-1)\big( |f'|^2\, +\, |f|^2 \big)}f'\ = \ f.
\eeq
Consequently, $f'$ is everywhere proportional to $f$ and as a result $f(\R)$ is contained into an $1$-dimensional subspace of $\R^N$. 

For $(ii)$, we have that if $k=1$ then \eqref{3.22} implies $f''=\la f + \mu f'$ for some $\la,\mu \in C^\infty(\R)$. Hence, $f(\R)$ is contained into a $2$-dimensional subspace of $\R^N$.  \eqref{3.21} gives the extra condition that $f'^\top(f'' +f)=0$ which implies $|f'|^2+|f|^2=c^2$  for some $c\in \R$. Hence, if $c\neq 0$ then $\frac{1}{c}(|f'|,|f|)^\top$ is on the unit circle and as such $|f|=c\cos B$ and $|f'|=c\sin B$, for some $B$ valued in $[0,\frac{\pi}{2}]$. Hence, $f=c\cos B (\cos A,\sin A)^\top$ for some $A$. The differential relation $|B'|^2 + |A'|^2 \cot^2 B =1$  follows easily.
\qed

\ms

\ms

\noi \textbf{Acknowledgement.} I am indebted to I.\ Vel\v{c}i\'c for our inspiring scientific discussions. I am also grateful to L.\ Capogna, Ch.\ Wang, J.\ Manfredi, F.\ Rindler, S.\ Aretakis and E.\ Scalas for their suggestions and constructive comments. I also thank L.\ C.\ Evans and Y.\ Yu for their interest and encourangement. Last but not least, I warmly thank the referees whose careful reading and constructive suggestions substantially improved both the content and the appearance of this paper.

\ms

\bibliographystyle{amsplain}

\end{document}